\documentstyle[amssymb,12pt]{article}
\begin{document}
\title{Polish group actions and admissible sets }
\author{B. Majcher-Iwanow}
\date{}
\maketitle

\setcounter{section}{-1} \newtheorem{thm}{Theorem}
\newtheorem{lem}[thm]
{Lemma} \newtheorem{definicja}[thm]{Definition}
\newtheorem{cor}[thm]
{Corollary} \newtheorem{prop}[thm]{Proposition}


\topskip 20pt

\begin{quote}
{\bf Abstract.}
We define some coding of Borel sets in admissible sets.
Using this we generalize certain results from model theory involving
admissible sets to the case of continuous actions of closed permutation 
groups on Polish spaces. 
In particular we obtain  counterparts of Nadel's theorems about
relationships between Scott sentences and admissible sets.

\bigskip



\end{quote}

\bigskip

\section{Introduction}

The aim of the paper is to study  actions of closed 
permutation groups on Polish spaces in admissible sets.
Let ${\mathbb{A}}$ be an admissible set.
Under some natural assumptions we can define in ${\mathbb{A}}$
a class of functions that can be considered as 'recipes' 
for Borel subsets of second countable spaces.
In Section 1 we describe such a coding and establish its basic properties.
\parskip0pt

Section 2 provides another tool of our study.
Let $G$ be a closed subgroup of $S_{\infty}$, the group
of all permutations of the set of natural numbers.
For every ordinal $\alpha<\omega_1$ we define $\alpha$-{\em sets},
Borel invariants that generalize on the one hand the concept of
a cannonical partition introduced by Becker in \cite{becker},
on the other - the concept of the $\alpha$-characteristic of 
a sequence in a structure given by Scott (see \cite{barwise}, p. 298).
In some other form these sets are defined and partially
studied by Hjorth in \cite{hjorth}.
Since they seem to be  interesting on its own rights
we examine their properties in detail.
We use them for the analysis of Borel complexity of $G$-orbits
(see Section 2.2).
\parskip0pt

The main results of the paper are proved in Section 3.
We prove that all Borel sets naturally involved in Scott analysis
can be coded in appropriate admissible sets.
Then we generalize Nadel's results concerning coding of
Scott analysis of countable structures in admissible sets \cite{nadel}.
We also give a generalization of another model
theoretical result - we characterize in admissible sets
orbits that are pieces of the canonical partition with respect to
some 'finer' topology ({\em nice topology} \cite{becker}).

A detailed description of our results is contained in Section 1.

\paragraph{Notation.}

A {\em Polish space (group)} is a separable, completely
metrizable topological space (group).
If a Polish group $G$ continuously acts on a Polish space $X$,
then we say that $X$ is a {\em Polish $G$-space}.
We usually assume that $G$ is considered under
a left-invariant metric.
We  say that a subset of $X$ is {\em invariant} if
it is $G$-invariant. \parskip0pt

We consider the group $S_{\infty}$ of all permutations
of the set $\omega$ of natural numbers and all its subgroups under the usual
left invariant metric $d$ defined by
$$
d(f,g)=2^{-\min\{k:f(k)\not=g(k)\}},\mbox{ whenever } f\not=g.
$$
We shall use the letters $a,b,c,d$ for finite sets of natural numbers.
For a finie set $d$ of natural numbers let $id_d$ be the identity map
$d\rightarrow d$ and $V_{d}$ be the group of all
permutations stabilizing $d$ pointwise, i.e.,
$V_d =\{f\in S_{\infty}: f(k)=k\mbox { for every } k\in d\}$.
Writing $id_{n}$ or $V_n$  we treat $n$ as the set of
all natural numbers less than $n$. \parskip0pt

Let $S_{<\infty}$ denote the set of all bijections between
finite substes of $\omega$.
We shall use small greek letters $\delta, \sigma, \tau$
to denote elements of $S_{<\infty}$.
For any  $\sigma \in S_{<\infty}$ let
$dom[\sigma], rng[\sigma]$ denote the domain
and the range of $\sigma$ respectively.

For every $\sigma \in S_{<\infty}$ let
$V_{\sigma} =\{f\in S_{\infty}:f\supseteq \sigma \}$.
Then for any $f\in  V_{\sigma}$
we have  $V_{\sigma}=fV_{dom[\sigma]}=V_{rng[\sigma]}f$.
Thus the family ${\mathcal N}=\{ V_{\sigma}:\sigma\in S_{<\infty}\}$
consists of all left (right) cosets of all subgroups $V_d$ as above.
This is a basis of the topology of $S_{\infty}$.

Given $\sigma\in S_{<\infty}$ and $s\subseteq dom[\sigma]$, then
for any $f\in V_{\sigma}$ we have $V_s^f=V_{\sigma [s]}$,
where $V_s^f$ denotes the conjugate $fV_sf^{-1}$.
\parskip0pt

In our paper we concentrate on Polish $G$-spaces,
where $G$ is a closed subgroup of $S_{\infty}$.
For such a group we shall use
the relativized version of the above, i.e.,
$V^G_{\sigma}=\{ f\in G:f\supseteq \sigma\}$,
$S_{<\infty}^G=\{ f|_d: f\in G$ and
$d$ is a finite set of natural numbers $\}$
(observe that for any subgroup $G$ and any finite set $d$
of natural numbers we have $id_{d}\in S_{<\infty}^G$).
The family
${\mathcal N}^G =\{ V_{\sigma}^G: \sigma\in S_{<\infty}^G\}$
is a basis of the standard topology of $G$. \parskip0pt

All basic facts concerning Polish $G$-spaces can
be found in \cite{bk}, \cite{hjorth} and \cite{kechris}.
\parskip0pt

Since we frequently use Vaught transforms,
recall the corresponding definitions.
The Vaught $*$-transform of a set $B\subseteq X$
with respect to an open $H\subseteq G$ is the set
$B^{*H}=\{ x\in X:\{ g\in H:gx\in B\}$ is comeagre in $H\}$,
the Vaught $\Delta$-transform of $B$ is the set
$B^{\Delta H}=\{ x\in X:\{ g\in H:gx\in B\}$
is not meagre in $H\}$.
It is known that for any $x\in X$ and $g\in G$, 
$gx\in B^{*H} \Leftrightarrow x\in B^{*Hg}$ and 
$gx\in B^{\Delta H} \Leftrightarrow x\in B^{\Delta Hg}$. 
On the other hand, if $B\in \Sigma^{0}_{\alpha}(X)$, 
then $B^{\Delta H}\in \Sigma^{0}_{\alpha}(X)$ 
and if $B\in\Pi^{0}_{\alpha}(X)$, then 
$B^{*H}\in\Sigma^{0}_{\alpha}(X)$.\parskip0pt  

It is worth noting that
{\em for any open $B\subseteq X$ and any open $K<G$ we have}
$B^{\Delta K}=KB$.
Indeed, by continuity of the action for any $x\in KB$ and
$g\in K$ with $gx\in B$ there are open neighbourhoods
$K_1\subseteq K$ and $B_1\subseteq KB$ of $g$ and $x$
respectively so that $K_1 B_1\subseteq B$; 
thus $x\in B^{\Delta K}$.
Other basic properties of Vaught transforms can be found
in \cite{bk} and \cite{kechris}. \parskip0pt

It is also assumed in the paper that the reader is
already acquainted with the most basic notions of
{\em admissible sets}.
Any necessary background can be easily provided by
\cite{barwise} and \cite{ershov}. \parskip0pt

We only remind the reader that an admissible set ${\mathbb{A}}$
is a transitive model of KPU, in the sense of \cite{barwise}.
Such models are considered as two-sorted structures
of some language $L$ with symbols $\emptyset ,\in$, where
one of the sorts corresponds to {\em urelements} and usually
forms a relational first-order structure with respect to the
symbols of $L$ distinct from $\emptyset$ and $\in$.
Here we assume that ${\mathbb{A}}$ satisfies KPU with respect
to all formulas of $L$ (${\mathbb{A}}$ is
{\em admissible with respect to} $L$ \cite{nadel}).


\section{Main results, Borel mulitcodes and Codability}

In this section we introduce the main notions of the paper 
and formulate our main results. 

To discuss Borel sets in an admissible set ${\mathbb{A}}$,
we shall assume that ${\mathbb{A}}$ contains some countable set 
(possibly as a set of urelements).
We will say that $\omega$ is {\em realizable} in an admissible set
${\mathbb{A}}$ if the set contains  a copy of the structure
$\langle \omega,<\rangle$ as an element (observe that $\omega$ 
is realizable in any admissible set satisfying Infinity Axiom).
If $\omega$ is realizable in an admissible set ${\mathbb{A}}$, then
by $\Delta$-separation ${\mathbb{A}}$ contains also a copy of
the set $[\omega]^{<\omega}$ of all finite sets of natural numbers,
a copy of $S_{<\infty}$ and, since $\subseteq$ is a $\Delta_0$-predicate,
copies of the posets $\langle [\omega]^{<\omega}, \subseteq \rangle$
and $\langle S_{<\infty}, \subseteq \rangle$.
Since it does not cause any misunderstanding, we shall write
$\omega$ and $S_{<\infty}$ even if we work not with the sets themselves 
but with their copies.

We start with the definition of Borel multicodes, i.e.
the functions that can serve as receipes for Borel sets.
Borel multicodes are not uniquely assigned to Borel sets, 
although every Borel multicode (with respect to a countable ordinal)
uniquely defines some Borel set. 

\begin{definicja} \label{D1}
Let ${\mathbb{A}}$ be an admissible set such that $\omega$ is realizable in it. 
We define in ${\mathbb{A}}$ two binary predicates $B_{\Sigma}$ and 
$B_{\Pi}$ by simultaneous induction on the ordinal $\alpha >0 $.
We put 
$$
\begin{array}{l@{\ }l}
B_{\Sigma} (1,u) & \mbox{ iff }\ u\mbox{ is a function }
\wedge dom[u]=\omega \wedge rng[u]\subseteq \{0,1\};\\

B_{\Pi} (\alpha,u) & \mbox{ iff }\  u=(0,u')\wedge
B_{\Sigma}(\alpha ,u');\\

\alpha>1\wedge B_{\Sigma} (\alpha ,u) & \mbox{ iff }\  u
\mbox{ is a function }
\wedge  (\alpha \mbox{ is a limit ordinal }\Rightarrow dom[u] = \alpha )\
\wedge \\

&\ \wedge \ (\alpha \mbox{ is the successor ordinal }\Rightarrow
dom[u]=\omega)\ \wedge\\

& \ \wedge \  (\forall u'\in rng[u])
        (\exists \beta <\alpha )(B_{\Pi}(\beta,u')\vee B_{\Sigma}(\beta,u'))
\end{array}
$$
If $\alpha$ is a non-zero ordinal then every $u$ such that
$B_{\Sigma}(\alpha, u)$ is called an $\alpha$-multicode
while every $u$ such that $B_{\Pi}(\alpha, u)$ is called
a co-$\alpha$-multicode.
\end{definicja}

We use some standard tricks of the general theory of definability 
in admissible sets (see \cite{barwise}) to show that
the relations above are $\Sigma$-definable.
Consider the  ternary predicate
$$
B(c,\alpha,u)\ \mbox{ iff }\  (c=0\wedge B_{\Sigma}(\alpha, u))\ \vee
\ (c=1\wedge B_{\Pi}(\alpha, u)).
$$
We see that the predicate $B(c,\alpha,u)$ is defined in
${\mathbb A}$ by a $B$-positive $\Sigma$-formula.
Thus by the second recursion theorem (Section 5.2 of \cite{barwise})
$B$ is a $\Sigma$-relation definable in ${\mathbb A}$.
Since $B_{\Sigma}(\alpha,u)$ is equivalent to $B(0,\alpha,u)$ and
$B_{\Pi}(\alpha,u)$ is equivalent to $B(1,\alpha,u)$, the predicates
$B_{\Sigma}$ and $B_{\Pi}$ are also $\Sigma$-predicates definable
in ${\mathbb A}$.

Now let ${\mathbb{A}}$ be an admissible set
such that $\omega$ is realizable in it.
Let $X$ be an arbitrary second countable space and 
$\{A_i:i\in \omega\}$ be its basis.
To every $u$ such that for some countable ordinal $\alpha \in{\mathbb{A}}$
we have ${\mathbb{A}}\models B_{\Sigma}(\alpha ,u)\vee B_{\Pi}(\alpha ,u)$,
we assign a Borel subset $B_u$ of $X$ in the following manner:
$$
\begin{array}{l@{\ \mbox{ then }\ }l}

\mbox{if}\quad  B_{\Sigma }(1,u) &  B_u=\bigcup\{A_n:u(n)=1\};\\

\mbox{if}\quad  B_{\Pi }(\alpha,u) & B_u=X\setminus B_{u'},
\mbox{ where }u=(0,u');\\

\mbox{if}\quad \alpha>1\wedge  B_{\Sigma} (\alpha,u) &
B_{u}=\bigcup \{B_{u'}: u'\in rng[u]\}.
\end{array}
$$ 
The assignment sends Borel multicodes $u$ satisfying
$B_{\Sigma}(\alpha,u)$ to the class $\Sigma_{\alpha}^0(X)$.
It is not one-to-one, in particular $B_u =B_v$, whenever
$B_{\Sigma}(\alpha ,u)$, $B_{\Sigma}(\alpha ,v)$ and $rng[u]=rng[v]$.

\begin{definicja} \label{D2}
Let ${\mathbb{A}}$ be an admissible set.
Let $X$ be a second countable space with a basis $\{ A_i:i\in \omega\}$
and $B\subseteq X$ be a Borel set.
If there are $u\in {\mathbb{A}}$ and a countable ordinal
$\alpha\in Ord ({\mathbb{A}})$ such that
${\mathbb{A}}\models B_{\Sigma}(\alpha,u)(\mbox{ or }
{\mathbb{A}}\models B_{\Pi}(\alpha, u))$ and $B=B_u$,
then we say that $B$ is constructible in ${\mathbb{A}}$ by $u$.
\end{definicja}

Observe that  the empty set, the whole space $X$ and
every basic open set $A_l$, are constructible by $1$-multicodes
in any admissible set ${\mathbb{A}}$ {\em realizing} $\omega$.
The functions $mc_{\emptyset}, mc_X, mc_l : \omega\to \{0,1\}$
below are the corresponding $1$-multicodes
$$
mc_{\emptyset}=(0,0,0,\ldots );\quad  mc_X=(1,1,1,\ldots );\quad
mc_l= (\underbrace{0,0,\ldots ,0}_{(l-1)-\mbox{times}},1,0,0,\ldots).
$$
We will use this notation below. \parskip0pt

Lemma \ref{constr} contains the most obvious properties
of {\em constructibility}.
In particular it states that this notion is preserved under
some natural operations which we shall use below.
Appropriate descriptions are given in the following definition.
By the second recursion theorem the predicate $Q_{\vee}$
defined below is a $\Sigma$-predicate.

\begin{definicja} \label{D3} 
Let ${\mathbb{A}}$ be an admissible set such that $\omega$ is realizable
in ${\mathbb{A}}$.
We define in ${\mathbb{A}}$ a ternary
predicate $Q_{\vee}$ by the following formula.
$$
Q_{\vee}(u,w,v )\Leftrightarrow (Q_0\wedge Q_1\wedge Q_2 )(u,w,v)
\mbox{ where }
$$
$$
\begin{array}{l@{\ }l}

Q_0 (u,w,v) & =  u, w, v\mbox{ are functions};\\

Q_1 (u,w,v) & = (\exists \alpha)(\alpha \mbox{ is a limit ordinal}
\wedge dom[u]=\alpha
\wedge dom[w]=\alpha \wedge dom[v ]=\alpha);\\

Q_2 (u,w,v) & =(\forall \beta<\alpha)(\forall n\in \omega )
\Big((\beta=0\vee \beta\mbox{ is a limit ordinal })\Rightarrow \\

&\Rightarrow (v(\beta+2n)=u(\beta +n)\wedge
v(\beta +2n+1)=w(\beta +n))\Big).
\end{array}
$$
\end{definicja}

It is easy to see that the predicate $Q_{\vee}$ defines an operation
on the class of all pairs of functions with common domain a limit ordinal.
We shall also use the following notation.
For any $u,w, v$ such that $Q_{\vee}(u,w,v)$ we shall write
$\bigvee (u,w)=v$.
If $u'=(0,u)$, $w'=(0,w)$ then we put $\bigwedge (u',w')=(0,\bigvee(u,w))$.
\parskip0pt

It is worth noting that if $\alpha$ is an ordinal and $u,w$ are
$\alpha$-multicodes then $\bigvee(u,w)$ is also an $\alpha$-multicode.
If $u,w$ are co-$\alpha$-multicodes then $\bigwedge (u,w)$ is also
a co-$\alpha$-multicode.

\begin{lem} \label{constr}
Let ${\mathbb{A}}$ be an admissible set and
$\alpha, \beta\in Ord({\mathbb{A}})$.
Let $X$ be a second countable space with a basis $\{ A_i :i\in \omega\}$
and $B, C\subseteq X$ be  Borel sets.\parskip0pt

(1) If  $\alpha <\beta$ and $B$ is constructible in
${\mathbb{A}}$ by some $u\in {\mathbb{A}}$ such that
${\mathbb{A}}\models B_{\Sigma}(\alpha, u)$ or
${\mathbb{A}}\models B_{\Pi}(\alpha, u)$ then there are
$w,w'\in {\mathbb{A}}$ such that

${\mathbb{A}}\models B_{\Sigma}(\beta, w)$
and ${\mathbb{A}} \models B_{\Pi}(\beta, w')$ and $B=B_w=B_{w'}$.
\parskip0pt

(2) If $B$ and $C$ are constructible in ${\mathbb{A}}$ by
some $\alpha$-multicodes $u$ and $w$ respectively then $B\cup C$
is constructible in ${\mathbb{A}}$ by $\bigvee (u,w)$;\parskip0pt

(3) If $B$ and $C$ are constructible in ${\mathbb{A}}$ by some
co-$\alpha$-multicodes $u$ and $w$ respectively then $B\cap C$
is constructible in ${\mathbb{A}}$ by $\bigwedge (u,w)$.
\end{lem}

{\em Proof.}
Let $u\in {\mathbb A}$ be an $\alpha$-multicode or a co-$\alpha$-multicode.
Then the function $w$ defined
by $w(n)=u$, for every $n\in \omega$
($w(\zeta )=u$, for every $\zeta < \beta$) is
a $\beta$-multicode for every successor (resp. limit)
ordinal $\beta>\alpha$.\parskip0pt

For turning $u$ into co-multicodes, note that the function
$u'$ defined by $u' (n)=(0,u)$ for all $n \in \omega$
($u'(\zeta )=(0,u)$ for every $\zeta< \beta$) satisfies
$B_{\Sigma }(\beta, z)$ and serves as a $\beta$-multicode
for $B_{(0,u)}$ for every  successor
(resp. limit) ordinal $\beta>\alpha$.
Then $w'$ can be taken as $(0,u')$. \parskip0pt

The rest of the lemma is easy.
$\Box$

\bigskip 

We now define some equivalence relation $\equiv$ on 
the set of multicodes (co-multicodes).

\begin{definicja} \label{eq}
Let ${\mathbb{A}}$ be an admissible set such that $\omega$ is realizable it.
We define in ${\mathbb A}$ a relation $\equiv$ by 
induction on the ordinal $\alpha >0 $: 
$$
u\equiv v\quad \mbox{ iff } \quad  \exists \alpha
\Big[\ (B_{\Sigma}(\alpha,u)\wedge B_{\Sigma}(\alpha,v))\ \wedge
(\ \alpha=1 \Rightarrow u=v\  )\ \wedge
$$
$$
(\alpha>1 \Rightarrow 
(\forall u'\in rng[u])(\exists v'\in rng[v])(u'\equiv v') \wedge 
(\forall v'\in rng[v])(\exists u'\in rng[u])(u'\equiv v')) 
$$
$$
\quad \vee\ \Big( B_{\Pi}(\alpha,u)\wedge B_{\Pi}(\alpha,v))
\wedge 2^{nd}(u)\equiv 2^{nd}(v)\Big)\ \Big]
$$
\end{definicja}

Since the operations $2^{nd}$, taking the second coordinate,
and $rng$, taking the range, are $\Sigma$-definable
(see Section 1.5 \cite{barwise}), we see that $\equiv$
is defined by a $\equiv$-positive $\Sigma$-formula.
Thus by the second recursion theorem it is
a $\Sigma$-relation in ${\mathbb{A}}$.
It is clear that $u\equiv v$ implies $B_u=B_v$.
The converse implication can fail.
On the other hand in some situations we will
be able to obtain some kind of this converse.
We will use it in Section 3 in
the proof of our main results. 

\bigskip

Now we are almost ready to discuss $G$-actions in admissible sets.
We only have to define some coding of information about an action
in admissible sets.

\begin{definicja}\label{F1}
Let $G<S_{\infty}$ be a closed subgroup and $\langle X, \tau\rangle$
be a Polish $G$-space with a basis $\{ A_l :l\in \omega\}$.
Let ${\mathbb{A}}$ be an admissible set.
We say that  $x\in X$ is {\em codable} (with respect to $G$)
in ${\mathbb{A}}$ if
$\omega$ is realizable in ${\mathbb{A}}$ and the function 
$$
F_1 : S_{<\infty}\to {\mathbb{A}}\quad\mbox{ defined by }\quad
F_1(\sigma)=\left\{\begin{array}{l@{\quad\mbox{if}\quad}l}
\emptyset&\sigma\not\in S_{<\infty}^G\\
\{l:V^G_{\sigma}  x\cap A_l \not=\emptyset\}&\sigma\in S_{<\infty}^G
\end{array}
\right.
$$
is an element of ${\mathbb{A}}$.
\end{definicja}

This condition corresponds to the standard assumption of 
\cite{nadel} that $M\in {\mathbb{A}}$ where $M$ is an element 
of the $S_{\infty}$-space of $L$-structures in the case of 
the logic action of $S_{\infty}$. 
In Section 3.2 we give a general straightforward construction which 
assigns an admissible set ${\mathbb{A}}_x$ to any element $x\in X$ 
such that $x$ is codable in ${\mathbb{A}}$. 

{\bf Remark.}
It is worth noting that in the definition we can demand only
that $F_1$ is $\Sigma$-definable in ${\mathbb{A}}$; 
then $F_1$ is an element of ${\mathbb{A}}$ by $\Sigma$-replacement 
(Theorem 1.4.6 from \cite{barwise}).
Using $\Delta$-separation (see \cite{barwise}, Theorems 1.4.5) 
we see that if $x$ is codable in ${\mathbb{A}}$ 
then the set
$S_{<\infty}^G=\{\sigma:\sigma\in S_{<\infty}, F_1(\sigma)\not=\emptyset\}$
is an element of ${\mathbb{A}}$. 
\parskip0pt 

In the situation when $x$ is codable in ${\mathbb{A}}$ 
we will usually assume that the relation
$$
Imp(c,l,k)\quad \Leftrightarrow \big( c\in [\omega ]^{<\omega}\wedge
l,k\in \omega\wedge \ A_k\subseteq V^G_c  A_l\big)
$$ 
is $\Sigma$-definable in ${\mathbb{A}}$.
This assumption is not very restrictive. 
For example when $X_L$ is the space of all $L$-structures 
on $\omega$ and $G=S_{\infty}$ acts on $X_L$ by
the logic action (see \cite{bk}), take any structure $M$ 
on $\omega$ with an appropriate coding of finite sets 
(for example the standard model of arithmetic).
Then ${\mathbb{A}} ={\mathbf{Hyp}}(M, Imp(c,l,k) )$, 
the admissible set above the structure $(M,Imp(c,l,k))$ 
has $Imp$  $\Delta_0$-definable (when $M=(\omega ,+,\cdot )$ we do
not even need to add $Imp$, because it is $\Sigma$-definable
in the structure).
In Section 3.2 we give some additional examples. \parskip0pt

The following theorem is the main result of the paper.

\begin{thm} \label{MaRe}
Let ${\mathbb{A}}$ be an admissible set such that
$\omega$ is realizable in it.
Let $G<S_{\infty}$ be a closed group, $X$ be a Polish $G$-space
with a basis $\{A_i:i>0\}$ and $Imp$ be $\Sigma$-definable on ${\mathbb{A}}$.
\parskip0pt

(1) Let $x\in X$ be $\Sigma$-codable in ${\mathbb{A}}$.
Then for every $y\in X$, if $x,y$ are in the same
invariant Borel subsets of $X$ which are constructible
in ${\mathbb{A}}$ then for every $\alpha\le o({\mathbb{A}})$
they are in the same invariant $\Sigma_{\alpha}^0$-subsets of $X$.
\parskip0pt

(2) If $x,y$ are $\Sigma$-codable in ${\mathbb{A}}$ and they
belong to the same invariant Borel sets which are constructible
in ${\mathbb{A}}$ then they are in the same $G$-orbit.
\end{thm}

It is based on Theorem \ref{con}, which will be proved in Section 3.
In fact the method is presented in Section 2, where
for every ordinal $\alpha<\omega_1$ we define $\alpha$-{\em sets}
$B_{\alpha}(x,\sigma )$, Borel invariants that generalize
on the one hand the concept of a canonical partition introduced by
Becker, on the other - the concept of an $\alpha$-characteristic
of a structure given by Scott.
$\alpha$-Sets appear in \cite{hjorth} in a slightly different form.
Since they seem to be interesting for its own rights we
examine their properties in detail.
Then we use them for the analysis of Borel complexity of $G$-orbits.
As a result we are able to improve several places of
Section 6.1 of \cite{hjorth}.
We also find a simplification of some theorem from
\cite{bk} on Borel orbit equivalence relations
in the case of actions of closed permutation groups.
\parskip0pt

It is worth noting that our results are not so
straightforward in the direction determined by Nadel.
Since we do not use standard tools from logic, we even
cannot formulate them in a sufficiently close form.
Instead of formulas (and of structures $\phi_{\alpha}$
used by Hjorth in \cite{hjorth}) we develope coding
of $\alpha$-sets $B_{\alpha}(x,\sigma )$ in admissible sets
(see Theorem \ref{con}).
As a result some fragments of Nadel's strategy look very
different in our approach.
In fact we completely avoid model theory in notation and
proofs. \parskip0pt

Theorem \ref{MaRe} suggests that under some additional
assumptions the orbit $Gx$ becomes the intersection
of all $G$-invariant Borel sets containing $x$ and
codable in ${\mathbb{A}}$.
In Section 3 we confirm this intuition
in the situation as follows.
Let  $(\langle X, \tau\rangle ,G)$ be
a Polish $G$-space with a countable basis
$\mathcal{A}$ consisting of clopen sets.
Along with the topology $\tau$ we shall consider
another topology on $X$.
The following definition comes from \cite{becker2}.

\begin{definicja} \label{nto}
A topology $t$ on $X$ is {\em nice} for the $G$-space
$(\langle X, \tau\rangle , G)$ if the following
conditions are satisfied.\\
(a) $t$ is a Polish topology, $t$ is finer than $\tau$
and the $G$-action remains continuous with respect to $t$.\\
(b) There exists a basis $\mathcal{B}$ for $t$ such that:\parskip0pt

(i) $\mathcal{B}$ is countable;\parskip0pt

(ii) for all $B_1, B_2 \in \mathcal{B}$, $B_1\cap B_2\in \mathcal{B}$;
\parskip0pt

(iii) for all $B\in \mathcal{B}$, $X\setminus B\in \mathcal{B}$;\parskip0pt

(iv) for all $B\in \mathcal{B}$ and $u\in \mathcal{N}^G$,
$B^{* u} \in \mathcal{B}$;\parskip0pt

(v) for any $B\in \mathcal{B}$ there exists an open subgroup $H < G$
such that $B$ is invariant \parskip0pt

under the corresponding $H$-action.\\
A basis satisfying condition $(b)$ is called a {\em nice} basis.
\end{definicja}

It is noticed in \cite{becker2} that any nice basis also
satisfies property (b)(iv) of the definition above
for $\Delta$-transforms.
It is also clear that
{\em any nice basis is invariant in the sense that
for every $g\in G$ and $B\in \mathcal{B}$ we have
$gB\in \mathcal{B}$} (see \cite{basia}). 
\parskip0pt

In Section 3 we will prove the following theorem.

\begin{thm}\label{mor}
Let $G$ be a closed subgroup of $S_{\infty}$,
$X$ be a Polish $G$-space, $t$ be a nice topology for $X$
and $\mathcal{B}$ be its nice basis.
Let $x\in X$ and let $C$ be the piece of the canonical partition
with respect to $\mathcal{B}$ containing $x$ (see \cite{becker}).
Let ${\mathbb{A}}$ be an admissible set such that
$x$ is codable in ${\mathbb{A}}$ with respect to $\mathcal{B}$.
Then the following are equivalent:\parskip0pt

(i) $C=G x$;\parskip0pt

(ii) $C$ can not be partitioned into two invariant Borel
sets constructible in ${\mathbb{A}}$.
\end{thm}

It is curious that this statement is related to some fact
from model theory, which was found by Morozov in \cite{morozov}.
Our proof is based on some arguments from \cite{basia}
together with the main tools of our paper.


\section{Sets arising in Polish group actions}

In this section we develope the generalized  Scott analysis
which was initiated in \cite{hjorth}.
We suggest a slightly different approach, more suitable
for the main tasks of the paper.
We replace the main tool of Hjorth's work
(hereditarily countable structures $\phi_{\alpha}(x,V_n )$ 
corresponding to Scott sentences)
by some invariants $B_{\alpha}(x,\sigma )$, $x\in X$,
$\sigma\in S^G_{<\infty}$, which are Borel subsets of the space.
They may be also used as counterparts of Scott sentences.
\parskip0pt 

Actually these sets already appear in \cite{hjorth}, where they are
defined in a different way 
\footnote{when $G=S_{\infty}$ it can be shown, that 
$B_{\alpha}(x,id_{n})=\{ y: \phi_{\alpha}(x,V_n )=\phi_{\alpha}(y,V_n )\}$
for $\alpha \ge\omega$} 
. 
We formulate another, more canonical definition.
It seems to be more convenient for many purposes.
It enables us to describe Borel complexity of the sets
$B_{\alpha}(x,\sigma )$ and compare it with Borel complexity
of the orbit $Gx$.
Finally they are more suitable for proofs of our main results
mentioned in the previous section.\parskip0pt

On the one hand this section can be considered as an improvement,
completion and systematization of the material scattered in 
Section 6.1 of \cite{hjorth}.
On the other hand it contains a couple of new results
(e.g. Propositions \ref{eff} and \ref{ef}) and a natural example,
which illustrates the introduced objects.

The section is divided into two subsections.
In the first one we define sets $B_{\alpha}(x,\sigma )$ and
describe the main properties of them.
Lemma \ref{min} is the key lemma which we use for
the main results of the paper.
On the other hand we study $\alpha$-sets $B_{\alpha}(x,\sigma )$
slightly further in order to present this material
in a complete form.
Propositions \ref{lemchar}, \ref{pol} and
\ref{ef} somehow summarize our study.
Proposition \ref{beke} (related to some results from \cite{bk})
is a straightforward application of our approach. \parskip0pt

In the second subsection we define a counterpart of the Scott
rank and compare it with the Borel rank of the orbit.

\subsection{Borel partitions}

Let $G$ be a closed subgroup of $S_{\infty}$ and $X$ be a Polish
$G$-space with a countable basis ${\mathcal{A}}=\{A_i:i\in \omega\}$.
We always assume throughout the paper
that every basic open set  is invariant with
respect to some basic clopen group $H<G$ (it follows from the continuity
of the action that such a basis exists).

By Proposition 2.C.2 of \cite{becker2} there exists a unique
partition of $X$, $X=\bigcup\{ Y_{t}: t\in T\}$ into invariant
$G_{\delta}$ sets $Y_{t}$ such that every orbit of $Y_{t}$ is
dense in $Y_{t}$.
To construct this partition we define for any $t\in 2^{\mathbb{N}}$ the set
$$
Y_{t}=(\bigcap\{ GA_{j}:t(j)=1\})\cap
(\bigcap\{ X\setminus GA_{j}:t(j)=0\})
$$
and take $T=\{ t\in 2^{\mathbb{N}}:Y_{t}\not=\emptyset\}$.\parskip0pt

In this section we generalize this notion and define for every ordinal
$0<\alpha<\omega_1$ some canonical partition of $X$  approximating
the original orbit partition.
In fact we define such  partitions not only for the whole group $G$,
but simultaneously for every basic clopen subgroup $V^G_d$, where $d$ is
a finite subset of $\omega$.
We call the classes of the partition
$\alpha$-{\em sets} and study their properties in detail.

\begin{definicja}\label{ph}
Let $G<S_{\infty}$ be a closed
subgroup and $\langle X, \tau\rangle$ be a Polish
$G$-space with a basis $\{ A_l :l\in \omega\}$.
For every $x\in X$  and $\sigma \in S_{<\infty}^G$ with $rng[\sigma]=c$
and $dom[\sigma]=d$ we define a Borel set $B_{\alpha}(x,\sigma)$ by
simultaneous induction on the ordinal $\alpha$.
$$
\begin{array}{l@{\ \ }l}
B_1(x,\sigma) & =\bigcap\{ V^G_{c} A_l:
V_{\sigma }^G x\cap A_l\not=\emptyset\}\cap
\bigcap\{X\setminus V^G_{c} A_l:
V_{\sigma}^G x\cap A_l=\emptyset\};\\

B_{\alpha +1}(x,\sigma)&=\bigcap\limits_{b\supseteq d} \
(\bigcup\{ B_{\alpha}(x,\sigma'):\sigma'\in S^G_{<\infty},
\sigma'\supseteq \sigma ,dom[\sigma' ]=b \})\cap\\
&\ \cap \bigcap\limits_{a\supseteq c}(\bigcup\{ B_{\alpha}(x,\sigma'):
\sigma'\in S^G_{<\infty},\sigma'\supseteq\sigma, rng[\sigma']=a\};\\

B_{\lambda}(x,\sigma) &=\bigcap\limits_{\alpha<\lambda}
B_{\alpha}(x,\sigma),\mbox{ for }\lambda \mbox{ limit } .
\end{array}
$$
\end{definicja}\bigskip

Although the definition of a $1$-set coincides with the definition
of a piece of the canonical partition, it is not quite evident that
the whole definition can be considered as a generalization of
the definition of the canonical partition.
This will be clearer when we describe some properties of the 
sets $B_{\alpha}(x, \sigma)$.
These properties will be applied in the main results of the paper.  

\begin{lem}\label{list}
Let $x,y\in X$, $\sigma\in S_{<\infty}^G$,
$dom[\sigma]=d$ and $rng[\sigma]=c$.
Then for any $f\in G$, $\delta \in S_{<\infty}^G$
and ordinals $\alpha ,\beta > 0$  the following statements are true.
\parskip0pt

(a) If $\beta\le\alpha$, then
$B_{\beta}(x,\sigma )\supseteq B_{\alpha}(x,\sigma)$;
\parskip0pt

(b) $B_{\alpha}(fx,\sigma)=B_{\alpha}(x, \sigma f)$, where
$\sigma f$ denotes the map $\sigma f|_{f^{-1}[d]}$, \parskip0pt

in particular  $B_{\alpha}(x,\sigma)=B_{\alpha}(fx,\sigma)$,
for every $f\in V^G_{d}$;\parskip0pt

(c) $fB_{\alpha}(x,\sigma)=B_{\alpha}(x,f\sigma)$,
in particular  $B_{\alpha}(x,\sigma)=B_{\alpha}(x,f\sigma)$,
for every $f\in V^G_{c}$;\parskip0pt

(d) $V^G_{\sigma} x \subseteq B_{\alpha}(x,\sigma)$ and
$B_{\alpha}(x,\sigma)$ is $V^G_{c}$-invariant;\parskip0pt

(e) $B_{\alpha +1}(x,\sigma)=\bigcap\limits_{\sigma' \supseteq \sigma}
V^G_{c} B_{\alpha}(x,\sigma')\ \cap
\ \bigcap\limits_{a\supseteq c}\bigcap\limits_{g\in V^G_c} \
(\bigcup\{ gB_{\alpha}(x,\sigma'):
\sigma'\in S^G_{\infty},\sigma'\supseteq\sigma ,$
\parskip0pt

$\ rng[\sigma']=a\});$
\parskip0pt

(f) If $\delta\supseteq \sigma$ then
$B_{\alpha}(x,\delta )\subseteq B_{\alpha}(x,\sigma)$;\parskip0pt

(g) If $y\in B_{\alpha}(x, \sigma)$  then
$B_{\alpha}(y,id_{c})=B_{\alpha}(x,\sigma)$;\parskip0pt

(h) If $rng[\delta]=c$  then either
$B_{\alpha}(x,\sigma)=B_{\alpha}(y,\delta)$ or
$B_{\alpha}(x,\sigma)\cap B_{\alpha}(y,\delta)=\emptyset$.
\end{lem}

{\em Proof.}
Statement of (a) follows directly from the definition.\parskip0pt

In the proof of (b) - (h) we shall frequently use
the following claim, which can be  derived by
easy straightforward arguments.\bigskip

{\em Claim}. Under the assumptions of the lemma we have:
\bigskip

1. $\{\{f\sigma':\sigma'\in S^G_{<\infty}, \sigma'\supseteq\sigma,
rng[\sigma']=a\}
: a\supseteq c\}=$\parskip0pt

$\ =\{\{\sigma':\sigma'\in S^G_{<\infty}, \sigma'\supseteq f\sigma,
rng[\sigma']=b\}: b\supseteq f[c]\}$; \bigskip

2. $\{\{f\sigma':\sigma'\in S^G_{<\infty}, \sigma'\supseteq\sigma,
dom[\sigma']=a\}
: a\supseteq d\}=$\parskip0pt

$\ =\{\{\sigma':\sigma'\in S^G_{<\infty}, \sigma'\supseteq
f\sigma, dom[\sigma']=b\}
: b\supseteq d\}$; \bigskip

3. If $f\in V^G_c$ then $\{\{f\sigma':\sigma'\in S^G_{<\infty},
\sigma'\supseteq\sigma, dom[\sigma']=a\}: a\supseteq d\}=$ \parskip0pt

$\ =\{\{\sigma':\sigma'\in S^G_{<\infty},\sigma'\supseteq \sigma,
dom[\sigma']
=b\}: b\supseteq d\}$;\bigskip

4. $\{\{\sigma'f:\sigma'\in S^G_{<\infty}, \sigma'\supseteq\sigma,
rng[\sigma']=a\}
: a\supseteq c\}=$\parskip0pt

$\ =\{\{\sigma':\sigma'\in S^G_{<\infty}, \sigma'\supseteq \sigma f,
rng[\sigma']=b\}: b\supseteq c\}$\bigskip

5. $\{\{\sigma'f:\sigma'\in S^G_{<\infty}, \sigma'\supseteq\sigma,
dom[\sigma']=a\}: a\supseteq d\}=$ \parskip0pt

$\ =\{\{\sigma':\sigma'\in S^G_{<\infty}, \sigma'\supseteq \sigma f,
dom[\sigma']=b\}: b\supseteq f^{-1}[d]\}$; \bigskip

6. If $f\in V^G_d$ then $\{\{\sigma'f:\sigma'\in S^G_{<\infty},
\sigma'\supseteq\sigma, dom[\sigma']=a\}: a\supseteq d\}=$ \parskip0pt 

$\ =\{\{\sigma':\sigma' \in S^G_{<\infty}, \sigma'\supseteq \sigma,
dom[\sigma']=b\}: b\supseteq d\}$;\bigskip

7. For each $\sigma' \in S^G_{<\infty}$ such that
$\sigma'\supseteq\sigma$ and $dom[\sigma']=b$  we have
$\{g\sigma':g\in V^G_c\}=\{\delta\in S^G_{<\infty}:
\delta\supseteq\sigma, dom[\delta]=b\}$.
\bigskip

Now we return to the proof of the lemma.\bigskip

(b) We proceed by induction on $\alpha>0$.
By the equality  $V^G_{\sigma f} x=V^G_{\sigma}f x$,
the stetement of (b) holds for $\alpha=1$.
Using the inductive assumption at the successor step  we get
$$
B_{\alpha+1}(fx,\sigma)=
$$
$$
\bigcap\limits_{b\supseteq d}\bigcup\{B_{\alpha}(x,\sigma'f):
\sigma' \in S^G_{<\infty},\sigma'\supseteq\sigma, dom[\sigma']=b\}\cap
$$
$$
\cap\bigcap\limits_{a\supseteq c}\bigcup\{B_{\alpha}(x,\sigma'f):
\sigma' \in S^G_{<\infty},\sigma'\supseteq\sigma,rng[\sigma']=a\}.
$$
Then we apply points 4 and 5 of the claim to get the required equality
$B_{\alpha+1}(fx,\sigma)=B_{\alpha+1}(x,\sigma f)$.
This completes the successor step.
The limit step is obvious.\bigskip

(c) By an obvious inductive argument we see that
$fB_{\alpha}(x,\sigma)=B_{\alpha}(fx,f\sigma f^{-1})$.
By (b) we obtain
$B_{\alpha}(fx,f\sigma f^{-1})=B_{\alpha}(x,f\sigma)$.
These equalities obviously imply the statement. \bigskip

(d) To prove the first part  we use induction on $\alpha$.
The inclusion trivially holds for $\alpha=1$.
The limit step is immediate.
Then we can easily settle the successor step,
since for every $b\supseteq d$ and $a\supseteq c$
we have
$$
V^G_{\sigma}=
\bigcup\{V^G_{\sigma'}:\sigma'\supseteq\sigma, dom[\sigma']=b\}
=\bigcup\{V^G_{\sigma'}:\sigma'\supseteq\sigma, rng[\sigma']=a\}.
$$\parskip0pt

The second part of (d) follows directly from (c).\bigskip

(e) By induction, using point (c) of the lemma and points 1, 3
of the claim.\bigskip

(f) We proceed inductively.
First we shall consider case $\alpha=1$.
If $A_l$ is a basic open set such that
$A_l\cap V^G_{\sigma} x\not=\emptyset$, then
$V^G_{\delta} x\subseteq V^G_{\sigma} x\subseteq V^G_c A_l$.
Since $V^G_c  A_l$ is open and the action is continuous,
there is a basic open set $A_k$ such that
$V^G_{\delta} x \cap A_k\not=\emptyset$ and
$V^G_{rng[\delta]} A_k\subseteq V^G_c A_l$.
This in particular implies that $B_1 (x,\delta)\subseteq V^G_c  A_l$.
\parskip0pt

On the other hand suppose that $A_l$ is a basic open set such that
$A_l\cap V^G_{\sigma} x=\emptyset$.
Since $V^G_{\sigma} x$ is $V^G_c$-invariant, we get
$V^G_c  A_l \cap V^G_{\sigma} x=\emptyset$.
We present $V^G_c  A_l$ as the union
$\ \bigcup\{ V^G_{rng[\delta ]} g A_l : g\in V^G_c\}\ $
and note that for every $g\in V^G_c$, we have
$V^G_{\delta} x \cap gA_l =\emptyset$.
Thus we have
$\ \bigcap \{ X\setminus V^G_{rng[\delta ]} A_k :V^G_{\delta}
x\cap A_k =\emptyset\}$
$\subseteq X\setminus V^G_c  A_l\ $ and then
$B_1 (x,\delta )\subseteq X\setminus V^G_c  A_l$.
This yields $B_1(x,\delta )\subseteq B_1(x,\sigma)$.\parskip0pt

For the successor step assume that the inclusion
$B_{\alpha}(x,\delta')\subseteq B_{\alpha}(x,\sigma')$ holds whenever
$\delta'\supseteq \sigma'$.
For any $b\supseteq d$ we put $\hat{b}=b\cup dom[\delta]$.
Using the inductive assumption we get
$$
\bigcup\{ B_{\alpha}(x,\delta'):\delta' \in S^G_{<\infty},
\delta'\supseteq\delta ,\ dom[\delta']=\hat{b}\}\subseteq
$$
$$
\subseteq\bigcup\{ B_{\alpha}(x,\sigma'):\sigma' \in S^G_{<\infty},
\sigma'\supseteq \sigma ,\ dom[\sigma']=\hat{b}\}\subseteq
$$
$$
\subseteq\bigcup\{ B_{\alpha}(x,\sigma'):\sigma' \in S^G_{<\infty},
\sigma'\supseteq\sigma ,\ dom[\sigma']=b\}.
$$
Similarly, if $a\supseteq c$ and $\hat{a}=a\cup rng[\delta]$ then
we have
$$
\bigcup\{ B_{\alpha}(x,\delta'):\delta' \in S^G_{<\infty},
\delta'\supseteq\delta ,\ rng[\delta']=\hat{a}\}\subseteq
$$
$$
\subseteq\bigcup\{ B_{\alpha}(x,\sigma'):\sigma' \in S^G_{<\infty},
\sigma'\supseteq\sigma ,\ rng[\sigma']=a\}.
$$
Hence we conclude that
$B_{\alpha +1}(x,\delta)\subseteq B_{\alpha +1}(x,\sigma)$.\parskip0pt

The limit step is immediate. \bigskip

(g) We proceed by induction.
For $\alpha =1$ the equality  follows directly from
the definition.
The limit step is immediate.
For the successor step, assume that the equality
$B_{\alpha}(x,\sigma')=B_{\alpha}(z,id_{a})$    holds
whenever $z\in B_{\alpha}(x,\sigma')$ and $rng[\sigma']=a$.
Now take an arbitrary
$y\in B_{\alpha +1}(x,\sigma)$.
By (e) we get
$$
B_{\alpha +1}(x,\sigma)=\bigcap\{V^G_c  B_{\alpha}(x,\sigma'):
\sigma'\in S^G_{<\infty} ,\ \sigma'\supseteq\sigma\} \cap
$$
$$
\cap \bigcap\limits_{a\supseteq c} \bigcap\limits_{g\in V^G_c}
(\bigcup\{ gB_{\alpha}(x,\sigma'):
\sigma'\in S^G_{<\infty} ,\ \sigma'\supseteq\sigma ,\ rng[\sigma']=a\}).
$$
We see that for every $\sigma'\in S^G_{<\infty}$ with
$\sigma'\supseteq \sigma$ there is some $f'\in V^G_{c}$
such that $f'y\in B_{\alpha}(x,\sigma')$ and thus
$B_{\alpha}(x,\sigma')=B_{\alpha}(y,id_{rng[\sigma']}f')$
(apply the inductive assumption and (c)).
Since $f'\in V^G_{c}$ and $\sigma'\supseteq \sigma$ then
$id_{rng[\sigma']}f'\supseteq id_{c}$ and
$rng[id_{rng[\sigma']}f']=rng[\sigma']$.
Hence the following is true
$$
(\forall \sigma'\supseteq \sigma)(\exists \delta'\supseteq id_{c})
\big( B_{\alpha}(x,\sigma')=B_{\alpha}(y,\delta') \wedge
rng[\sigma']=rng[\delta']\big) .
$$
On the other hand take an arbitrary $\delta'\supseteq id_{c}$.
Put $a=rng[\delta']$ and take any $g\in V^G_c$ such that
$g\supseteq \delta'$.
Then by (f) there is some $\sigma'\supseteq \sigma$ such that
$rng[\sigma']=a$ and $gy\in B_{\alpha}(x,\sigma' )$.
By the inductive assumption, the latter implies
$B_{\alpha}(x,\sigma')=B_{\alpha}(gy,id_{a})$.
Then by (c) we get $B_{\alpha}(x,\sigma')=B_{\alpha}(y,id_{a} g)
=B_{\alpha}(y,\delta')$.
\parskip0pt

We have proved that if $y\in B_{\alpha +1}(x,\sigma)$
then the equality
$$
\{ B_{\alpha}(x,\sigma'):\sigma' \in S^G_{<\infty},
\sigma'\supseteq \sigma ,\ rng[\sigma']=a\}=
\{ B_{\alpha}(y,\delta'):\delta' \in S^G_{<\infty},\delta'\supseteq id_{c} ,\ rng[\delta']=a\}
$$
is true for every $a\supseteq c$.
This implies
$$
\bigcap\limits_{a\supseteq c}(\bigcup\{ B_{\alpha}(x,\sigma'):
\sigma' \in S^G_{<\infty},\sigma'\supseteq\sigma, rng[\sigma']=a\}=
$$
$$
=\bigcap\limits_{a\supseteq c}
(\bigcup\{ B_{\alpha}(y,\delta'):
\delta' \in S^G_{<\infty},\delta'\supseteq id_c, rng[\sigma']=a\}
$$ and
$$
\bigcap\{V^G_c B_{\alpha}(x,\sigma'):\sigma' \in S^G_{<\infty},
\sigma'\supseteq \sigma\}=
\bigcap\{V^G_c B_{\alpha}(y,\delta'):
\delta' \in S^G_{<\infty},\delta'\supseteq id_c\}
$$
which by (e) gives the required equality.\bigskip

(h) follows directly from (g). $\Box$
\bigskip

We are now ready to prove that the partition of $X$ 
into $\alpha$-sets can be defined by the same scheme 
as the canonical partition (thus can be considered as 
a generalization of the latter). 

\begin{prop}\label{canon} 
Let ${\mathcal A}$ be a basis for $X$,
$x\in X$,  $\sigma\in S^G_{<\infty}$, $rng[\sigma]=c$ 
and $\alpha>0$ be an ordinal.\parskip0pt

(a) Let ${\cal B}_{\alpha}=
\{B_{\alpha}(y,\delta):y\in X, \delta\in S^G_{<\infty}\}$.
Then we have
$$
B_{\alpha+1}(x,\sigma)=
\bigcap\{V^G_c B:
B\in {\mathcal B}_{\alpha},\ V^G_{\sigma} x\ \cap \ B\not=\emptyset\}\cap
$$
$$
\cap\bigcap\{X\setminus  V^G_c B:
\ B\in {\mathcal B}_{\alpha},\ V^G_{\sigma} x\ \cap \ B=\emptyset\}.
$$
\parskip0pt

(b) Let ${\mathcal B}_{<\alpha}=
\{B_{\gamma}(y,\delta):y\in X, \delta\in S^G_{<\infty}, \gamma<\alpha\}
\cup {\mathcal A}.$
Then we have
$$
B_{\alpha}(x,\sigma)=
\bigcap\{V^G_c B:
B\in {\mathcal B}_{<\alpha},\ V^G_{\sigma} x\ \cap \ B\not=\emptyset\}\cap
$$
$$
\cap\bigcap \{X\setminus  V^G_c B:
\ B\in {\mathcal B}_{<\alpha},\ V^G_{\sigma} x\ \cap \ B=\emptyset\}.
$$
\end{prop}

{\em Proof}. (a) The inclusion $\supseteq$ easily follows from the definition
and the lemma above.
We have to work a little more with its  converse.
Let $B\in {\mathcal B}_{\alpha}$ be such that
$V^G_{\sigma} x\ \cap \ B\not=\emptyset$.
Then , by the lemma above,
there is some $\sigma' \in S^G_{<\infty}$ such that
$\sigma'\supseteq \sigma$ and $B_{\alpha}(x,\sigma')\subseteq B$.
Hence  we have
$V^G_c B_{\alpha}(x,\sigma')\subseteq V^G_c B$, which yields
$B_{\alpha+1}(x,\sigma)\subseteq V^G_c B$.\parskip0pt

On the other hand let $B\in {\mathcal B}_{\alpha}$ be such that
$V^G_{\sigma} x\ \cap \ B=\emptyset$.
Then , by the lemma above,
there is some $a\supseteq c$ such that
$B_{\alpha}(x,\sigma')\cap B=\emptyset$, for every
$\sigma' \in S^G_{<\infty}$ with $\sigma'\supseteq \sigma$ and
$rng[\sigma']=a$.
Therefore
$$
\bigcap\limits_{a\supseteq c}\bigcup \{B_{\alpha}(x,\sigma'):
\sigma' \in S^G_{<\infty},\sigma'\supseteq \sigma, rng[\sigma']=a\}
\cap V^G_c B=\emptyset ,
$$
which yields
$B_{\alpha+1}(x,\sigma)\subseteq X\setminus V^G_c B$.
\parskip0pt

(b) follows from (a)  and the properties of $\alpha$-sets
collected in Lemma \ref{list}. $\Box$
\bigskip

Proposition \ref{canon} (b) yields the folowing statement.

\begin{prop}\label{=}
Let $x,y\in X$, $\alpha>1$ be an ordinal and $c\subseteq \omega$ be
a finite set.
Then for every $\sigma,\delta\in S^G_{\infty}$ with common range $c$
the following are equivalent:
\parskip0pt

(i) $B_{\alpha}(x,\sigma)=B_{\alpha}(y,\delta)$;\parskip0pt

(ii) For every finite $a\supseteq c$ and every $\zeta<\alpha$ we have
\parskip0pt

$\quad \{B_{\zeta}(x,\sigma'):\sigma'\supseteq\sigma, rng[\sigma']=a\}=
\{B_{\zeta}(y,\delta'):\delta'\supseteq\delta, rng[\delta']=a\};$\parskip0pt

(iii) For every natural $n\supseteq c$ and every $\zeta<\alpha$ we have
\parskip0pt

$\quad \{B_{\zeta}(x,\sigma'):\sigma'\supseteq\sigma, rng[\sigma']=n\}=
\{B_{\zeta}(y,\delta'):\delta'\supseteq\delta, rng[\delta']=n\}.$
\end{prop}

{\em Proof.} By Proposition \ref{canon}, (i) is equivalent to the equality
$$
\{B\in {\mathcal B}_{<\alpha}: B\cap V^G_{\sigma} x\not=\emptyset\}=
\{B\in {\mathcal B}_{<\alpha}: B\cap V^G_{\delta} y\not=\emptyset\}.
$$
(i) $\Rightarrow$ (ii)
Fix arbitrary $\zeta<\alpha$ and $a\supseteq c$.
Then take any $\sigma'\supseteq\sigma$ with $rng[\sigma']=a$.
Since $B_{\zeta}(x,\sigma')\cap V^G_{\sigma}x\not=\emptyset$,
we see that
$B_{\zeta}(x,\sigma')\cap V^G_{\delta}y\not=\emptyset$.
Hence by Lemma \ref{list} (g), there is some
$\delta'\supseteq \delta$ with $rng[\delta']=a$ such that
$B_{\zeta}(x,\sigma')=B_{\zeta}(y,\delta')$.
Therefore
$
\{B_{\zeta}(x,\sigma'): \sigma' \in S^G_{<\infty},
\sigma'\supseteq\sigma, rng[\sigma']=a\}\subseteq
\{B_{\zeta}(y,\delta'): \delta' \in S^G_{<\infty}, rng[\delta']=a\}.
$
In the same way we derive the converse inclusion.
\parskip2pt

(ii)$\Rightarrow $(i)
Take any $B\in {\mathcal B}_{<\alpha}$ such that
$B\cap V^G_{\sigma}x\not=\emptyset$.
There are $\zeta<\alpha$ and $\sigma'\subseteq \sigma$ such that
$B_{\zeta}(x,\sigma')\subseteq B$.
Since we can find $\delta'\supseteq\delta$ (with $rng[\delta']=rng[\delta]$)
such that $B_{\zeta}(x,\sigma')=B_{\zeta}(y,\delta')$,
we see that  $B\cap V^G_{\delta}y\not=\emptyset$.
This proves
$
\{B\in {\mathcal B}_{<\alpha}: B\cap V^G_{\sigma} x\not=\emptyset\}
\subseteq \{B\in {\mathcal B}_{<\alpha}: B\cap V^G_{\delta} y\not=\emptyset\}
$.
Similarly we obtain the converse inclusion.
\parskip2pt

(ii) $\Rightarrow$ (iii) is obvious.\parskip2pt

To prove (iii) $\Rightarrow$ (ii) suppose that (ii) does not hold.
Then there are some finite set $a\supseteq c$ and ordinal $\zeta<\alpha$
such that
$$
\{B_{\zeta}(x,\sigma'):\sigma' \in S^G_{<\infty},
\sigma'\supseteq\sigma, rng[\sigma']=a\}\not=
\{B_{\zeta}(y,\delta'):\delta' \in S^G_{<\infty},
\delta'\supseteq\delta, rng[\delta']=a\}.
$$
Take any natural $n\supseteq a$.
By Lemma \ref{list} (f), (h), we have
$$
\{B_{\zeta}(x,\sigma'):\sigma'\supseteq\sigma, rng[\sigma']=n\}\not=
\{B_{\zeta}(y,\delta'):\delta'\supseteq\delta, rng[\delta']=n\},
$$
hence (iii) does not hold.
$\Box$
\bigskip

The lemma below shall play the key role in the proof of
the main result of the paper. 
It states that $\alpha$-sets are in some sense minimal 
with respect to $\alpha$. 

\begin{lem}\label{min}
Let $x\in X$, $\sigma\in S_{<\infty}^G$,
$c=rng[\sigma]$ and $\alpha>0$ be an ordinal.
Then for any $V^G_c$-invariant
$A\in \Sigma^0_{\alpha}\cup\Pi^0_{\alpha}$\quad we have
\quad  $V^G_{\sigma} x\subseteq A\quad \mbox{iff}\quad
B_{\alpha}(x,\sigma)\subseteq A.$
\end{lem}

{\em Proof.}
To prove ($\Rightarrow$) we proceed inductively.\parskip0pt

Consider case $\alpha=1$.
If $U$ is a $V^G_c$-invariant open set containing
$V^G_{\sigma} x,$ then there is a basic open set $A_{l_0}\subseteq U$
intersecting $V^G_{\sigma} x.$
Then $V^G_c\, A_{l_0}\subseteq U$ and so
$$
U\supseteq \bigcap\{ V^G_c A_l :A_l\cap V^G_{\sigma} x\not=\emptyset\}
\supseteq B_1(x,\sigma).
$$
If $F$ is an $V^G_c$-invariant closed set containing
$V^G_{\sigma} x$, then
$$
F\supseteq\bigcap \{X\setminus V^G_c\, A_l: F\cap A_l=\emptyset\}\supseteq
\bigcap \{X\setminus V^G_c\, A_l: V^G_{\sigma} x\cap A_l=\emptyset\}
\supseteq B_1(x,\sigma).
$$
The rest of the proof is based on the following statements.
\parskip4pt

{\em Claim } Let $\alpha$ be an ordinal, $c\in [\omega]^{<\omega}$
and $A\subseteq X$ be a $V^G_c$-invariant set.
\parskip3pt

(1) If $A\in \Sigma^0_{\alpha}(X)$
then  $A$ can be presented as a union $A=\bigcup\limits_i D_i$ such that
$\{D_i:\ i<\omega\}\subseteq \bigcup\limits_{\xi<\alpha}\Pi^0_{\xi}(X)$
and for every $i<\omega$ there is $a_i\supseteq c$  such that
$D_i$ is  $V^G_{a_i}$-invariant. Moreover, if $\alpha$ is limit then
each $D_i$, $\ i<\omega\ $, can be taken $V^G_c$-invariant.\parskip3pt

(2) If $A\in \Pi^0_{\alpha}(X)$ is a $V^G_c$-invariant set
then  $A$ can be presented as an intersection
$A=\bigcap\limits_i D_i$ such that
$\{D_i:\ i<\omega\}\subseteq \bigcup\limits_{\xi<\alpha}\Sigma^0_{\xi}(X)$
and for every $i<\omega$ there is $a_i\supseteq c$  such that
$D_i$ is  $V^G_{a_i}$-invariant. Moreover, if $\alpha$ is limit then
each $D_i$, $\ i<\omega\ $, can be taken $V^G_c$-invariant.\parskip3pt

{\em Proof of Claim.}
(1) There is a countable family
$\{A_i:i\in \omega\}\subseteq \bigcup\limits_{\xi<\alpha}\Pi^0_{\xi}(X)$
such that
$A=\bigcup\limits_i A_i$.
Since $A$ is $V^G_c$-invariant we have
$$
A=A^{\Delta V^G_c}=\bigcup\limits_i
(\bigcup\{ A^{* W}_i :W\subseteq V^G_c \mbox{ is basic, open}\} )=
$$
$$
=\bigcup\limits_i \bigcup\limits_{a\supseteq c}
(\bigcup\{ A^{* V^G_{\delta}}_i :(\delta\in S^G_{<\infty})
\wedge (\delta\supseteq id_{c})\wedge (dom[\delta]=a)\} ).
$$
It follows from the properties of Vaught transforms that if
$A_i\in \Pi^0_{\xi}$ and $dom[\delta]=a$ then
$A^{* V^G_{\delta}}_i$ is a $V^G_a$-invariant
$\Pi^0_{\xi}$-set.
It completes the first part.\parskip1pt

Now, it is clear that if $A_i\in \Pi^0_{\xi}$ 
then the set   
$$
A_i^{\Delta V^G_c}=\bigcup\limits_{a\supseteq c}
(\bigcup\{ A^{* V^G_{\delta}}_i :(\delta\in S^G_{<\infty})
\wedge (\delta\supseteq id_{c})\wedge (dom[\delta]=a)\} ) 
$$ 
is a $V^G_c$-invariant $\Sigma_{\xi+1}$-set.
Thus it is also a $V^G_c$-invariant $\Pi_{\xi+2}$-set.
Then $A$ is a countable union of $V^G_c$-invariant  elements of the union
$\bigcup\limits_{\xi<\alpha}\Pi^0_{\xi+2}(X)$, which proves the additional
statement for limit $\alpha$.\parskip2pt

(2) There is a countable family
$\{A_i:i\in \omega\}\subseteq \bigcup\limits_{\xi<\alpha}\Sigma^0_{\xi}(X)$
such that $A=\bigcap\limits_i  A_i$.
Since $A$ is $V^G_c$-invariant, we have
$$
A=A^{* V^G_c}=\bigcap\limits_i (\bigcap\{ A_i^{\Delta W}:
W\subseteq V^G_c\mbox{ is basic, open}\} )=
$$
$$
=\bigcap\limits_i \bigcap\limits_{a\supseteq c}
(\bigcap\{ A_i^{\Delta V^G_{\delta}} : (\delta \in S^G_{<\infty})
\wedge (\delta\supseteq id_{c})\wedge (dom[\delta ]=a)\} ).
$$
Applying standard properties of Vaught transform again, we see that
if $A_i\in \Sigma^0_{\xi}$ and $dom[\delta]=a$.
Hence $A^{\Delta V^G_{\delta}}_i$ is a $V^G_a$-invariant
$\Sigma^0_{\xi}$-set.
It completes the first part.\parskip1pt

Now,  if $A_i\in \Sigma^0_{\xi}$ then the set
$$
A_i^{*V^G_c}=\bigcap\limits_{a\supseteq c}
(\bigcap\{ A^{* V^G_{\delta}}_i :(\delta\in S^G_{<\infty})
\wedge (\delta\supseteq id_{c})\wedge (dom[\delta]=a)\} ) 
$$ 
is a $V^G_c$-invariant $\Pi_{\xi+1}$-set,
thus it is also a $V^G_c$-invariant $\Sigma_{\xi+2}$-set.
Then $A$ is a countable intersection of $V^G_c$-invariant  elements of
the union $\bigcup\limits_{\xi<\alpha}\Sigma^0_{\xi+2}(X)$,
which proves the additional
statement for limit $\alpha$.\bigskip

We continue the proof of Lemma.

To go through the successor step
take an arbitrary  $\alpha$ and assume that the statement
holds for every $a\in [\omega]^{<\omega}$, $\sigma'\in S_{<\infty}^G$
with $a=rng[\sigma']$ and every  $V^G_{a}$-invariant set
$D\in \Sigma^0_{\alpha}\cup\Pi^0_{\alpha}$
containing $V^G_{\sigma'} x$.
Let $A\in \Sigma^0_{\alpha+1}(X)\cup \Pi^0_{\alpha+1}(X)$
be an arbitrary $V^G_c$-invariant set containing $V^G_{\sigma} x$.
We shall consider two cases.\bigskip

$1^o$ $A\in \Sigma^0_{\alpha+1}$.\parskip0pt

By Claim, $A$ can be presented as a union
$A=\bigcup\limits_i D_i$, where
for every $i<\omega$ there is $a_i\supseteq c$  such that
$D_i$ is  a $V^G_{a_i}$-invariant $\Pi^0_{\alpha}$-set.\parskip0pt

Fix an arbitrary $g\in V^G_{\sigma}$.
Since $V^G_{\sigma} x\subseteq A$,
there are $i\in \omega$ and $a_i\supseteq c$ such that $gx\in D_i$ and
$D_i$ is $V^G_{a_i}$-invariant.
Put $\sigma' =id_{a_i}g$.
Then we have $\sigma'\supseteq \sigma$, $rng[\sigma']=a_i$ and
$V^G_{\sigma'} x\subseteq  D_i$.
Using the inductive assumption we conclude that
$B_{\alpha}(x,\sigma')\subseteq D_i\subseteq A$.
Since $A$ is $V^G_c$-invariant we obtain
$B_{\alpha+1}(x,\sigma)\subseteq V^G_c  B_{\alpha}(x,\sigma' )
\subseteq A$.\bigskip

$2^o$ $A\in \Pi^0_{\alpha+1}$.
By Claim, $A$ can be presented as a union
$A=\bigcap\limits_i D_i$, where
for every $i<\omega$ there is $a_i\supseteq c$  such that
$D_i$ is  a $V^G_{a_i}$-invariant $\Sigma^0_{\alpha}$-set.\parskip0pt

Fix arbitrary $i\in \omega$ and $a_i\supseteq c$ such that
$D_i$ is $V^G_{a_i}$-invariant.
We have
$$
\bigcup\{ V^G_{\sigma'} x :(\sigma'\supseteq \sigma )\wedge (rng[\sigma']
=a_i)\}= V^G_{\sigma} x\subseteq D_i.
$$
Thus for every $\sigma'\in S^G_{<\infty}$ with
$\sigma'\supseteq \sigma$ and $rng[\sigma']=a_i$ we have
$V^G_{\sigma'} x \subseteq D_i$.
Since $D_i$ is a $V^G_{a_i}$-invariant
$\Sigma^0_{\alpha}$-set, by the inductive assumption we conclude
that $B_{\alpha}(x,\sigma')\subseteq D_i$.
Therefore
$$
\bigcup\{ B_{\alpha}(x,\sigma'): (\sigma'\in S^G_{<\infty})\wedge
(\sigma'\supseteq\sigma)\wedge (rng[\sigma']=a)\}
\subseteq D_i.
$$
By  Definition \ref{ph}  this completes the successor step.\parskip0pt

By Claim we can also easily go through the limit step.

The backward direction is just Lemma \ref{list} (d).
$\Box$
\bigskip

The next result provides  another  necessary and sufficient condition
for the equality of $\alpha$-sets.
It improves the result from \cite{hjorth}, where some counterpart
of (i) $\Rightarrow$ (ii) is proved.

\begin{prop} \label{lemchar}  
Let $x,y\in X$, $\sigma,\delta\in S_{<\infty}^G$ and
$rng[\sigma]=rng[\delta]=c$.
Then for every ordinal $\alpha>0$ the following conditions are equivalent.
\parskip0pt

\quad(i) $B_{\alpha}(x,\sigma)=B_{\alpha}(y,\delta);$\parskip0pt

\quad (ii) For every $V^G_c$-invariant set
$A\in \Sigma^0_{\alpha}(X)\cup \Pi^0_{\alpha}(X)$
we have\parskip0pt

\quad\quad
$\ V^G_{\sigma} x\subseteq A\ \mbox{ iff }
\ V^G_{\delta} y\subseteq A.$

Moreover for every ordinal $\alpha>0$ we have
$
B_{\alpha}(x,\sigma)=\bigcap\{A\in \Sigma^0_{\alpha}(X)\cup
\Pi^0_{\alpha}(X):\ A\mbox { is }
V^G_c\mbox{-invariant},\  V^G_{\sigma} x\subseteq A\}.$

In particular $B_{\alpha}(x,\sigma)$ is a $\Pi^0_{\alpha +1}$-set
for every successor ordinal $\alpha$, and  a $\Pi^0_{\alpha}$-set
for every limit ordinal $\alpha$. \parskip0pt
\end{prop}

{\em Proof.}
$(i)\Rightarrow (ii)$ follows from Lemma \ref{list}(d)
and Lemma \ref{min}.\parskip0pt

To prove $(ii)\Rightarrow (i)$ we use induction on $\alpha>0$.
The case $\alpha=1$ and the limit step are easy.
To go through the successor step, take an arbitrary $\alpha$ and assume
that for  every $\sigma', \delta'$
with $rng[\sigma']=rng[\delta']=a$ if
$B_{\alpha}(x,\sigma')\not=B_{\alpha}(y,\delta')$ then
we can separate $V^G_{\sigma'} x$ from $V^G_{\delta'} y$ by some
$V^G_a$-invariant set $A\in \Sigma^0_{\alpha}\cup \Pi^0_{\alpha}$.
Then suppose that
$B_{\alpha+1}(x,\sigma)\not=B_{\alpha +1}(y,\delta)$.
By Proposition \ref{=}  there is some $a\supseteq c$
such that one of the following cases holds
\begin{quote}
$1^o$
For some $\sigma'\supseteq \sigma$ with $rng[\sigma']=a$ 
and every $\delta'\supseteq \delta$  with $rng[\delta']=a$ we have
$B_{\alpha}(x,\sigma')\not=B_{\alpha}(y,\delta')$;\parskip0pt

$2^o$ For some $\delta'\supseteq \delta$ with $rng[\delta']=a$
and every $\sigma'\supseteq \sigma$  with $rng[\sigma']=a$ we have
$B_{\alpha}(x,\sigma')\not=B_{\alpha}(y,\delta')$.
\end{quote}
Since the cases are symmetric we consider only the first one.
By the inductive assumption, for every
$\delta'\supseteq \delta$  with $rng[\delta']=a$
there is some $V^G_{a}$-invariant set
$A_{\delta'}\in \Sigma^0_{\alpha}(X)\cup \Pi^0_{\alpha}(X)$
such that $V^G_{\sigma'} x\subseteq A_{\delta'}$
while  $V^G_{\delta'} y$ is disjoint from $A_{\delta'}$.
Then for every $\sigma'\supseteq\sigma$ with $rng[\sigma']=a$ we have
$$
V^G_{\sigma'} x\subseteq \bigcap\{ A_{\delta'}:
(\delta' \in S^G_{<\infty})\wedge (\delta'\supseteq \delta)\wedge
(rng[\delta']=a)\}
$$
while the set
$V^G_{\delta} y=
\bigcup\{ V^G_{\delta'}y:\delta'\supseteq \delta\wedge rng[\delta']=a\}$
is disjoint from
$\bigcap\{ A_{\delta'}: \delta'\supseteq \delta
\wedge rng[\delta']=a\}$.

Put
$A=\bigcap\{ A_{\delta'}^{* V^G_c}: \delta'\supseteq \delta
\wedge rng[\delta']=a\}$.
Then $A\in \Pi^0_{\alpha +1}(X)$,
$V^G_{\sigma} x = (V^G_{\sigma'} x)^{* V^G_c}\subseteq  A$
and $V^G_{\delta} y$ is disjoint from  $A$.

The second part of the statement is a direct consequence
of the previous one and Lemma \ref{min}. $\Box$ 
\bigskip 

As an immediate consequence of this proposition and Lemma \ref{list}(h)
we obtain the following statement.

\begin{cor} \label{complex}
Let  $\sigma\in S_{<\infty}^G$ and $rng[\sigma]=c$.
\parskip0pt

(a) For every  successor ordinal $\alpha$ the family
$\{B_{\alpha}(x,\sigma):x\in X\}$
is a partition of $X$ into $V^G_c$-invariant $\Pi^0_{\alpha+1}$-sets.
\parskip0pt

(b) For every limit ordinal $\alpha$ the family
$\{B_{\alpha }(x,\sigma):x\in X\}$ is a partition of $X$ into
$V^G_c$-invariant $\Pi^0_{\alpha}$-sets.
\end{cor}

\bigskip

Every piece of the canonical partition as a $G_{\delta}$-subset of $X$,
is a Polish space with the topology inhertited from the original
Polish topology on $X$.
We  generalize this fact and show that every $\alpha$-set is a Polish space
with respect to some finer topology generated by 'ealier' $\beta$-sets.
From now on we shall use the following notation for every ordinal $\beta$:
$$
\begin{array}{l}
{\cal B}^x_0={\cal A}\\
{\cal B}^x_{\beta}= \{B_{\beta}(x,\sigma): \sigma\in S^G_{<\infty}\}
\mbox{\quad for\quad }\beta>0\\
{\cal B}^x_{<\beta}=\bigcup\{{\cal B}^x_{\gamma}:\ \gamma<\beta \}.
 \end{array}
$$

\begin{prop}\label{pol}
Let ${\cal A}$ be a countable basis of $X$,
$x\in X$  and $0<\alpha<\omega_1$ be an ordinal.
The set $B_{\alpha}(x,\emptyset)$
with the (relative) topology generated by the family
${\cal B}^x_{<\alpha}$ as basic open sets is a Polish $G$-space.
\end{prop}

{\em Proof.}
As we have already mentioned, $B_1(x,\emptyset)$ is a $G_{\delta}$ subset
of $X$, thus is a Polish space with respect to the (relative) topology 
generated by ${\cal A}$.
Therefore we will  deal below only with $\alpha>1$.
We shall use the following result by Sami (see \cite{sami}, Lemma 4.2).
\begin{quote}
{\em Let $\langle X, t\rangle$ be a topological space and
$1\le \zeta<\omega_1$.
Let ${\cal F}$ be a {\em Borel family of rank} $\zeta$, i.e.
a family of subsets of $X$ which can be decomposed into
subfamilies of two types
$\ {\cal F}=\bigcup\{P_{\xi}:1\le\xi<\zeta\}\cup
\bigcup\{S_{\xi}:1\le\xi<\zeta\}
\ $
satisfying the following conditions:
1. $S_1$ consists of open sets,
2. $P_{\xi}=\{X\setminus A:\ A\in S_{\xi}\}$, for $1\le \xi<\zeta$,
3. every element of $S_{\xi}$ is a union of a countable subfamily
of $\bigcup \{P_{\eta}:1\le\eta<\xi\}$, for $1\le \xi<\zeta$.\parskip0pt

If $X$ is a Polish space then the topology generated by a family
of intersections of finite subsets of the union $t\cup {\cal F}$ is also
Polish.}
\end{quote}

Consider the family
$$
{\hat {\cal B}}^x_{<\alpha}=
\left\{\begin{array}{l}
\{ V_c B,\  X\setminus V_c B:\ B\in {\cal B}^x_{<\beta},
\  c\in [\omega]^{<\omega}\}\cup
\{ B, X\setminus B: B\in {\cal B}^x_{\beta}\}\\
\mbox{ if }\ \alpha=\beta+1>1\\
\\
\{ V_c B,\  X\setminus V_cB:\ B\in {\cal B}^x_{<\alpha},
\  c\in [\omega]^{<\omega}\}\mbox{ if }\ \alpha\mbox{ is a limit ordinal }.
\end{array}\right.
$$
Our first task is to show that $\hat{\cal B}^x_{<\alpha}$ is 
a Borel family of some countable rank.
To prove this we need some preliminary work. 

We define for every $1<\xi<\omega_1$ the sets $S_{\xi}$ and $P_{\xi}$.
First we put:
$$
\begin{array}{l}
S_1=\{V_c A: A\in {\cal A}, c\in [\omega]^{<\omega}\}\\
P_1=\{X\setminus D:\  D\in S_1\}\\
\\
S_2=P_1\\
P_2=\{X\setminus D: \ D\in S_2\}\\
\\
S_3=\Big\{\ \bigcup\{ V_c A:A\cap V^G_{\sigma}x=\emptyset\}\cup
\bigcup\{X\setminus V_c A:A\cap V^G_{\sigma}\not=\emptyset\}: \\
\ \ \ \ \ c\in [\omega]^{<\omega}, \sigma\in S^G_{<\infty}, rng[\sigma]=c\  \Big\}\\
P_3=\{X\setminus D:\ D\in S_3\}
\end{array}
$$
Observe that
$\bigcup\limits_{i=1}\limits^{3}(S_i\cup P_i)$
is a Borel family of rank 4.

We proceed similarly at each successor stage.
Every  successor ordinal has one of
following form:   $\xi + 3n+1$, $\xi +3n+2$ or $\xi +3n+3$,
where $n$ is a natural number and
$\xi=0$ or $\xi$ is a limit ordinal.
We define
\bigskip

$
\begin{array}{l}
S_{\xi +3n+1}=\{ V_c B: B\in {\cal B}^x_{\xi+n}\}\\
P_{\xi +3n+1}=
\{X\setminus V_c B: B\in {\cal B}^x_{\xi +n}\}\\
\\
S_{\xi +3n+2}=P_{\xi +3n+1}\\
P_{\xi +3n+2}=S_{\xi +3n+1}\\
\\
S_{\xi +3n+3}=\{X\setminus B:\ B\in {\cal B}^x_{\xi+n+1}\}\\
P_{\xi+3n+3}={\cal B}^x_{\xi+n+1}
\end{array}
$
\bigskip

Finally, for every limit $\xi<\omega_1$ we put:\bigskip

$
\begin{array}{l}
S_{\xi}=\{X\setminus B:\ B\in {\cal B}^x_{\xi}\}\\
P_{\xi}={\cal B}^x_{\xi}.
\end{array}
$\bigskip

We claim that for every $1<\zeta<\omega_1$ the family
$\bigcup\limits_{\xi<\zeta}(S_{\xi}\cup P_{\xi})$ is a Borel family
of rank $\zeta$.
It is clear that such a family satisfies conditions 1-2.
We have to check that  it also satisfies condition 3.
We apply an inductive argument.
It is obvious for $\zeta=2, 3, 4$.
The case of a limit $\zeta$ is immediate either.
For the successor step take an arbitrary $\zeta>1$ and
suppose that the family
$\bigcup\limits_{0<\xi<\zeta}(S_{\xi}\cup P_{\xi})$
satisfies condition 3.
We have  $\bigcup\limits_{0<\xi<\zeta+1}(S_{\xi}\cup P_{\xi})=
\bigcup\limits_{0<\xi<\zeta}(S_{\xi}\cup P_{\xi})
\ \cup \ S_{\zeta}\cup P_{\zeta}$.

Consider two cases.\parskip0pt

$1^o$ $\zeta$ is limit.
By Definition \ref{ph} we have
$$
\begin{array}{l}
P_{\zeta}={\cal B}^x_{\zeta}
=\Big\{ \bigcap\limits_{\xi<\zeta} B_{\xi}(x,\sigma):
\sigma\in S^G_{<\infty}\Big\}\\
S_{\zeta}=\{X\setminus B:\ B\in {\cal B}^x_{\zeta}\}
=\Big\{ \bigcup\limits_{\xi<\zeta} (X\setminus B_{\xi}(x,\sigma)):
\sigma\in S^G_{<\infty}\Big\} .
\end{array}
$$
By the definition of the sets $S_{\xi}, P_{\xi}$
and the assumption that $\zeta$ is limit  we see that
${\cal B}^x_{<\zeta}\subseteq\bigcup\limits_{\xi<\zeta}P_{\xi}\subseteq
\bigcup\limits_{\xi<\zeta} S_{\xi}$.
Hence $\{ X\setminus B:\  B\in {\cal B}^x_{<\zeta }\}
\subseteq\{X\setminus B:\  B\in \bigcup\limits_{\xi<\zeta}S_{\xi}\}\subseteq
\bigcup\limits_{\xi<\zeta}P_{\xi}$.
Therefore every element of $S_{\zeta}$ is a countable union of elements
of the set $\bigcup\limits_{\xi<\zeta}P_{\xi}$ which completes Case $1^o$.
\parskip2pt

$2^o$ $\zeta$ is a successor ordinal.

There are unique ordinals $\gamma$ and $n$ such that
$n$ is a  natural number,
$\gamma$ equals $0$ or is a limit ordinal and
$\zeta$  has one of the following form:
$\gamma+3n+1$, $\gamma+3n+2$ or $\gamma+3n+3$.
If $\zeta$ takes one of the first two forms, then we are done
directly from the definition.\parskip0pt

If $\zeta=\gamma+3n+3$ then by Proposition \ref{canon} we have
$$P_{\zeta}={\cal B }^x_{\gamma+n+1}=
\Big\{\bigcap\{V_c B:\ B\in {\cal B}^x_{\gamma+n},
B\cap V^G_{\sigma}x\not=\emptyset\}\cap$$
$$
\cap \bigcap\{X\setminus V_cB:\ B\in {\cal B}^x_{\gamma+n},
B\cap V^G_{\sigma}=\emptyset\}:\\
\ c\in [\omega]^{<\omega}, \sigma\in S^G_{<\infty}, rng[\sigma]=c\Big
\}
$$
and
$$
S_{\zeta}=\{X\setminus B: \ B\in {\cal B}^x_{\gamma+n+1}\}
=\Big\{\bigcup\{V_c B:\ B\in {\cal B}^x_{\gamma+n},
B\cap V^G_{\sigma}=\emptyset\}\cup
$$

$$
\cup\bigcup\{X\setminus V_cB:\ B\in {\cal B}^x_{\gamma+n},
B\cap V^G_{\sigma}\not=\emptyset\}:\\
\ c\in [\omega]^{<\omega}, \sigma\in S^G_{<\infty}, rng[\sigma]=c\Big
\}.
$$
Since $\{V_c B:\ B\in {\cal B}^x_{\gamma+n}\}=P_{\gamma+3n+2}$
and $\{X\setminus V_cB:\ B\in {\cal B}^x_{\gamma+n}\}=P_{\gamma+3n+1}$
we conclude that $S_{\zeta}$ consists of countable unions
of elements from $P_{\zeta+3n+1}\cup P_{\zeta+3n+2}$.
This completes Case $2^o$.\bigskip

Now let $\gamma$ and $k$ be the unique ordinals such that
$\alpha=\gamma+k$, $k$ is a natural number and
$\gamma$ equals $0$ or is a limit ordinal.
Define
$$\hat{\alpha}=
\left\{\begin{array}{l@{\ \mbox{if}\ }l}
\gamma+3(k-1)+1 & k>0\\
\alpha & k=0.\end{array}\right.
$$
We see that
$\hat{\cal B}^x_{<\alpha}=\bigcup\{P_{\xi}:1\le\xi<\hat{\alpha}\}\cup
\bigcup\{S_{\xi}:1\le\xi<\hat{\alpha}\}$.
Hence ${\hat {\cal B}}^x_{<\alpha}$
is a Borel family of  rank $\hat{\alpha}$.\bigskip

Since $\hat{{\cal B}}^x_{<\alpha}$ is countable,
it generates a Polish topology on $X$.
Since $B_{\alpha}(x,\emptyset)$ is a $G_{\delta}$-subset of $X$
with respect to this topology, it is a Polish space with
the inherited topology.
As we have already noted
${\cal B}^x_{<\alpha}\subseteq {\hat{\cal B}}^x_{<\alpha}$.
We now show that every set of the form
$B_{\alpha}(x,\emptyset)\cap D$, where $D\in {\hat {\cal B}}^x_{<\alpha}$
is a union of elements from
$\{B_{\alpha}(x,\emptyset)\cap B:\ B\in {\cal B}^x_{<\alpha}\}$,
i.e. the latter family can be also taken as a basis of the topology.
It follows from the following claim.\bigskip

{\em Claim.} Let $\zeta<\beta<\alpha$.
Then for every $\sigma\in S^G_{<\infty}$ and $c\in [\omega]^{<\omega}$
the sets $B_{\alpha}(x,\emptyset)\cap V_cB_{\zeta}(x,\sigma)$
and  $B_{\alpha}(x,\emptyset)\setminus V_cB_{\zeta}(x,\sigma)$
are unions of elements from the family
$\{B_{\alpha}(x,\emptyset)\cap B:\ B\in  {\cal B}^x_{\beta}\}$.\parskip0pt

{\em Proof of Claim.} Take any
$y\in B_{\alpha}(x,\emptyset)\cap V_cB_{\zeta}(x,\sigma)$.
By Proposition \ref{canon} we get
$B_{\beta}(y,id_c)\subseteq V_cB_{\zeta}(x,\sigma)$.
On the other hand Lemma \ref{list}(h) yields
$B_{\alpha}(y,\emptyset)=B_{\alpha}(x,\emptyset)$.
Then, by Proposition \ref{=} we conclude that there is some
$B\in {\cal B}^x_{\beta}$ such that $B=B_{\beta}(y,id_c)$
which proves the first part of the claim.\parskip0pt

Similarly, if
$y\in B_{\alpha}(x,\emptyset)\setminus V_cB_{\zeta}(x,\sigma)$, then
on the one hand
$B_{\beta}(y,id_c)\subseteq X\setminus V_cB_{\zeta}(x,\sigma)$,
on the other hand
$B_{\alpha}(y,\emptyset)=B_{\alpha}(x,\emptyset)$.
Using  Proposition \ref{=} again we conclude that for some
$B\in {\cal B}^x_{\beta}$ we have
$y\in B\cap B_{\alpha}(x,\emptyset)\subseteq
B_{\alpha}(x,\emptyset)\setminus  V_cB_{\zeta}(x,\sigma)$,
which proves the second part.\bigskip

Now it suffices to notice that the action
$a: G \times B_{\alpha}(x,\emptyset)\to B_{\alpha}(x,\emptyset)$
is continuous with respect to each argument.
Since every element of ${\cal B}^x_{<\alpha}$ is
invariant with respect to a basic open subgroup of $G$,
the action is continuous with respect to the first coordinate.
On the other hand, for every $f\in G$, $\gamma<\alpha$ and
$\delta\in S^G_{<\infty}$  we have
$\{y\in X: fy\in B_{\gamma}(x,\delta)\}=B_{\gamma}(x,f^{-1}\delta)$,
which proves continuity with respect to the second coordinate.
$\Box$
\bigskip

From now on let $t^x_{\alpha}$ denote the Polish topology on
$B_{\alpha}(x,\emptyset)$ described above.
Observe that in the case when $\alpha$ is a successor ordinal
and $\alpha=\beta+1$, the topology $t_{\alpha}^x$ is also (relatively)
generated by a smaller basis, namely ${\cal B}^x_{\beta}$.
It follows directly from the claim used in the proof above and
Corollary \ref{complex}.\parskip3pt

Using Proposition \ref{pol} together with Effros Theorem on $G_{\delta}$-orbits 
we obtain the following fact. 

\begin{prop}\label{eff}
Let $G$ be a closed subgroup of $S_{\infty}$, $X$ be a Polish $G$-space 
and $x\in X$.
Let $\alpha >0$ be an ordinal.\parskip0pt

(a) If $G x\in \Pi^0_{\alpha+1}(X)$,
then the  following statements are true:\parskip0pt

\quad (i) $G x$ is non-meager in $B_{\alpha}(x,\emptyset)$ with respect
to $t^x_{\alpha}$.\parskip0pt

\quad (ii) the map\quad
$G\to G x\quad g\to g x $
\quad  is open with respect to $t^x_{\alpha}$.\parskip0pt

\quad (iii)  $G x= B_{\alpha+1}(x,\emptyset)$.\parskip0pt

(b) If $\alpha$ is a limit ordinal and  $G x\in \Pi^0_{\alpha}(X)$,
then the  following statements are true:\parskip0pt

\quad (iv)  $G x= B_{\alpha}(x,\emptyset)$.\parskip0pt

\quad (v) the map\quad
$G\to G x\quad g\to g x $
\quad  is open with respect to $t^x_{\alpha}$.
\end{prop}

{\em Proof}. (a) The proof is based on the following observation. 
\parskip3pt

{\em Claim. } Let $A\subseteq B_{\alpha}(x,\emptyset)$ be an invariant
$\Pi^0_{\alpha+1}$ -set.
Then $A$ is a $G_{\delta}$-set with respect to $t^x_{\alpha}$.

{\em Proof of Claim. } We apply the claim used in the proof of
Lemma \ref{min}.
We present $A$ as countable intersection  $A=\bigcap\limits_{i<\omega} D_i$,
such that for every $i<\omega$, $D_i\in \Sigma^0_{\alpha}(X)$
and it is invariant with respect $V^G_{a_i}$ for  some finite
$a_i\subseteq \omega$.
Then we apply the claim again to each $D_i$.
For every $i<\omega$, we find a family  $\{D_{ij}:\ j<\omega\}$ satisfying
the following conditions:\parskip2pt

1. $D_{ij}\in \bigcup\limits_{\xi<\alpha}\Pi^0_{\xi}(X)$;\parskip0pt

2. $D_{ij}$ is invariant with respect to $V^G_{a_{ij}}$, for
some finite $a_{ij}\supseteq a_i$;\parskip0pt

3. $D_i=\bigcup\limits_{j<\omega} D_{ij}$.\parskip2pt

We have
$$
A=A\cap B_{\alpha}(x,\emptyset)=
\bigcap\limits_{i<\omega}\bigcup\limits_{j<\omega}
(D_{ij}\cap B_{\alpha}(x,\emptyset)).
$$

By Lemma \ref{min} and Lemma \ref{list}(h) we see
that for every $i,j<\omega$ the set
$D_{ij}\cap B_{\alpha}(x,\emptyset)$ is  a union 
of elements of the family
$\{B\cap B_{\alpha}(x,\emptyset):\ B\in {\cal B}^x_{<\alpha}\}$,
thus it is open with respect to the topology $t^x_{\alpha}$.
Therefore $A$ is a $G_{\delta}$ set with respect to $t^x_{\alpha}$.
\parskip4pt

Hence we conclude that
$G x$ is a $G_{\delta}$-subset of $B_{\alpha}(x,\emptyset)$ 
with respect to the topology.
Then (i) and (ii) follows from Effros theorem.
(iii) follows  Lemma \ref{min}.\parskip3pt

(b) Point (iv) follows from Lemma \ref{min}, then we obtain  (v)  
from Effros theorem. $\Box$

\bigskip

The second statement of the following proposition looks folklore, 
but we have not found it in literature.

\begin{prop}\label{ef} Let $G$ be a closed subgroup of $S_{\infty}$,
$X$ be a Polish $G$-space.\parskip0pt

(a) Let $x\in X$. If $G x\in \Pi^0_{\alpha}(X)$
for some ordinal $\alpha$, then for every open basic subgroup $V^G_c<G$
we have $V^G_c x\in \Pi^0_{\alpha}(X)$,
and for every open subgroup
$H<G$ we have $H x\in \Pi^0_{\alpha+2}(X)$  .\parskip0pt

(b) The orbit equivalence relation induced on $X$ by the $G$-action
is Borel if and only if the orbit equivalence relation induced on $X$ 
by the action of some of open subgroup $H<G$ is Borel.
\end{prop}

{\em Proof}. (a) Let $V^G_c$ be an arbitrary basic open subgroup of G.
 We have to consider two cases.\parskip0pt

1$^o$ $\alpha=\beta+1$ is a successor ordinal.
Then by Proposition \ref{eff}(ii) there is a family
${\mathcal C}\subseteq {\mathcal B}^x_{<\beta}$ such that
$$
(*) \ V^G_c x = \big(\bigcup {\mathcal C}\big)\cap G x.
$$
Thus for some $C_0\in {\cal C}$ we have
$C_0\cap V^G_c x\not= \emptyset$, 
which implies $B_{\beta}(x,id_c)\subseteq V^G_c C_0$.
Next, since the set on the right side of the equality $ (\star)$
must be $V^G_c$-invariant, we have also
$V^G_c x = V^G_c\big( \bigcup {\mathcal C}\big)\cap G x$.
The latter implies 
$B_{\beta}(x,id_c)\cap G x\subseteq V^G_cC_0\cap G x \subseteq  V^G_c x$, 
which yields $V^G_c x=B_{\beta}(x,id_c)\cap G x$. 
Then we are done, since ${\mathcal B}^x_{\beta}\subseteq \Pi_{\alpha}(X)$.
\parskip0pt

2$^o$ $\alpha$ is a limit ordinal.
Then by Proposition \ref{eff}(v) there is a family
${\mathcal C}\subseteq {\mathcal B}^x_{<\alpha}$ such that
$$
(\star) \ V^G_c x = \big(\bigcup {\mathcal C}\big)\cap G x.
$$
Exactly as in the case $\alpha =\beta +1$ we 
obtain $V^G_c x=B_{\alpha}(x,id_c)\cap G x$.
Then we are done, since 
${\mathcal B}^x_{\alpha}\subseteq \Pi_{\alpha}(X)$
for every limit $\alpha$. \parskip0pt

The second part of (a) is a direct consequence of the first one.

(b) is a consequence of (a), the fact that every $G$-orbit 
is a countable union of $H$-orbits and the following
theorem of Sami on Borel orbit equivalence relations
(see \cite{sami}). 
\begin{quote}
Let $G$ be a Polish group and $X$ be a Polish $G$-space.
The orbit equivalence relation induced on $X$ by the $G$-action is 
Borel if and only if there is a countable ordinal $\alpha$ such
that every $G$-orbit is a $\Pi^0_{\alpha}$-subset of $X$. 
\end{quote}
$\Box$
\bigskip

We now show that Proposition \ref{ef} simplifies the proof 
of Theorem 7.1.1 of Becker and Kechris from \cite{bk} 
in the particular case of actions of closed subgroups of
$S_{\infty}$ (in \cite{bk} it is assumed that $G$ is Polish). 
In fact S.Solecki suggested that such applications are possible. 

\begin{prop} \label{beke} (Becker, Kechris - special case) 
Let $G$ be a closed subgroup of $S_{\infty}$ and 
$X$ be a Borel $G$-space.
Then the following conditions are equivalent.\parskip0pt

(i) The orbit equivalence relation is Borel.\parskip0pt

(ii) The map $\tau: x\to G_x$ from $X$ to
the Effros space of closed subsets ${\cal F}(G)$ is Borel.
\end{prop}

{\em Proof.}
Note that since the action is continuous, all stabilizers of $G$
are closed.
For (ii)$\Rightarrow $ (i) see \cite{bk}.
To prove the converse we can assume that $X$ is a Polish $G$-space.
Each basic open set in ${\cal F}(G)$ has the form
$U_{\sigma}=\{K\in {\cal F}(G): K\cap V^G_{\sigma}\not=\emptyset\}$,
where $\sigma\in S^G_{\infty}$.
Take an arbitrary  $\sigma\in S^G_{\infty}$  and fix some $g\in V^G_{\sigma}$.
Then $V^G_{\sigma}=V^G_cg$, where $c=rng[\sigma]$.
Let $E_c$ denote the orbit equivalence relation induced by $V^G_c$.
We  have $\tau^{-1}[U_{\sigma}]=\{x\in X: G_x\cap V^G_{\sigma}\not=\emptyset\}
=\{x\in X:(x,gx)\in E_c\}=\pi_X[\{(x,gx):x\in X\}\cap E_c]$.
By Proposition \ref{ef}, $E_c$ is Borel.
Then we are done since the projection $\pi_X$ is one-to-one on the set
$\{(x,gx):x\in X\}\cap E_c$. $\Box$

\subsection{Ranks of orbits} 

Now we shall define for every $x\in X$ some cardinal 
invariant connected with $\alpha$-sets. 
The definition is based on the following lemma. 

\begin{lem}\label{height}
For every $x\in X$ there is some $\gamma <\omega_1$
such that for all $\sigma, \delta \in S_{<\infty}^G$ with
$rng[\sigma]=rng[\delta]$  we have
$$
(\exists \alpha<\omega_1)\big(B_{\alpha}(x,\sigma)\not=B_{\alpha}(x,\delta )
\big)
\Rightarrow\big(B_{\gamma}(x,\sigma)\not=B_{\gamma}(x,\delta)\big).
$$
\end{lem}

{\em Proof.} Let $d\subseteq \omega$ be an arbitrary finite set.
For every pair $\{\sigma,\delta\}\subseteq S^G_{<\infty}$ with
the same range $d$ consider the set
$\{\alpha <\omega_1:B_{\alpha}(x,\sigma)\not=B_{\alpha}(x,\delta)\}$.
Let $\gamma_{\sigma,\delta}$ be its infimum in case the set is nonempty or
$0$ otherwise.
Let $\gamma_d=sup\{\gamma_{\sigma,\delta}:
\{\sigma,\delta\}\subseteq S^G_{<\infty}, rng[\sigma]=rng[\delta]=d\}$.
It is a countable ordinal, since the set
$\{\{\sigma,\delta\}:rng[\sigma]=rng[\delta]=d\}$ is countable.
Finally let $\gamma=sup\{\gamma_d:d\in [\omega]^{<\omega}\}$.
It is also a countable ordinal as a supremum of
a countable set of countable ordinals.
Obviously the ordinal $\gamma$ has the required property.
$\Box$

\begin{definicja}
For every $x\in X$ let $\gamma^G_{\star}(x)$ be the least ordinal
$\gamma$ satisfying the statement of Lemma \ref{height}.
\end{definicja}

By Lemmas \ref{min} and \ref{list}  every orbit is an $\alpha$-set.
In the theorem below we show that every $G x$ is a
$(\gamma^G_{\star}(x)+2)$-set.

\begin{thm}\label{orb}
For every $x\in X$ we have
$B_{\gamma^G_{\star}(x)+2}(x,\emptyset)=G x$.
\end{thm}

{\em Proof.} Let $y\in B_{\gamma^G_{\star}(x)+2}(x,\emptyset)$.
Then by Lemma \ref{list} (h) we have
$B_{\gamma^G_{\star}(x)+2}(x,\emptyset)=
B_{\gamma^G_{\star}(x)+2}(y,\emptyset)$.
The rest of the proof is based on two claims.\bigskip

{\em Claim 1}. Let $\zeta, \beta$ be ordinals such that
$\gamma^G_{\star}(x)< \zeta+1<\beta$ and
$B_{\beta}(x,\emptyset)=B_{\beta}(y,\emptyset)$.
Let $\sigma,\delta\in S^G_{<\infty}$ have common range $a$.
Then the equality $B_{\zeta}(x,\sigma)=B_{\zeta}(y,\delta)$  implies
$B_{\zeta+1}(x,\sigma)=B_{\zeta+1}(y,\delta)$.\parskip0pt

{\em Proof of Claim 1}. By Proposition \ref{=} there is some $\sigma'$
such that $B_{\zeta+1}(x,\sigma')=B_{\zeta+1}(y,\delta)$.
Then by Lemma \ref{list}(a), (i) we have also
$B_{\zeta}(x,\sigma')=B_{\zeta}(y,\delta)$  and
$B_{\zeta}(x,\sigma)=B_{\zeta}(x,\sigma')$.
Since $\zeta\ge\gamma^G_{\star}$,
we have $B_{\zeta+1}(x,\sigma)=B_{\zeta+1}(x,\sigma')$,
which proves the required equality.
\bigskip

{\em Claim 2}. If $B_{\gamma^G_{\star}(x)+2}(x,\emptyset)=
B_{\gamma^G_{\star}(x)+2}(y,\emptyset)$ then
$B_{\beta}(x,\emptyset)=B_{\beta}(y,\emptyset)$
for every $\beta$.\parskip0pt

{\em Proof of Claim 2}.
It is trivially true for $\beta\le \gamma^G_{\star}+2$.
We use induction to prove it for ordinals  $\beta>\gamma^G_{\star}+2$.
\parskip0pt

The limit step is immediate.
To go through the successor step suppose
that $B_{\beta}(x,\emptyset)=B_{\beta}(y,\emptyset)$, for some
$\beta\ge \gamma^G_{\star}+2$.
Take an arbitrary $a\subseteq \omega$.
Then by Proposition \ref{=}, for every $\zeta<\beta$ we have
$\{B_{\zeta}(x,\sigma): rng[\sigma]=a\}=\{B_{\zeta}(y,\delta):
rng[\delta]=a\}$.
This equality remains true for $\zeta=\beta$.
Indeed, it is obvious for a limit $\beta$
and follows from Claim 1 for a successor $\beta$.
By Proposition \ref{=} again we get
$B_{\beta+1}(x,\emptyset)=B_{\beta+1}(y,\emptyset)$. \bigskip

We come back to the proof of the theorem.
From the fact that $G$-orbits are Borel sets, we conclude
by Lemma \ref{min}, that there is an ordinal $\beta$ such that
$G x=B_{\beta}(x,\emptyset)$ and $G y=B_{\beta}(y,\emptyset)$.
Since 
$B_{\gamma^G_{\star}(x)+2}(x,\emptyset)=B_{\gamma^G_{\star}(x)+2}(y,\emptyset)$, 
we are done by Claim 2. 
$\Box$ 
\bigskip

As a corollary of the theorem we obtain the following statement.

\begin{cor}
For every $x\in X$ and $\sigma\in S^G_{<\infty}$ we have
$B_{\gamma^G_{\star}(x)+2}(x,\sigma)=V^G_{\sigma} x$.
\end{cor}

{\em Proof.} 
We derive it from the definition of $\gamma^x_{\star}$ and
Lemma \ref{list} (d), (f). $\Box$\bigskip

The following lemma  gives a characterization of
$\gamma^G_{\star}$ in terms of Borel complexity.
It is a direct consequence of Proposition \ref{lemchar} 
(together with the idea that appropriate Vaught 
transforms make a set invariant). 

\begin{lem}\label{char}
For every $x\in X$, $\gamma^G_{\star}(x)$ is the least ordinal $\alpha$ with
the property that
\begin{quote}
for every $\sigma, \sigma_1\in S^G_{<\infty}$ with
$\ rng[\sigma]=rng[\sigma_1]\ $ if
$V^G_{\sigma} x \not=V^G_{\sigma_1} x$ then there is
a Borel set of rank $\alpha$ containing one of the set
$V^G_{\sigma} x$, $V^G_{\sigma_1} x$
and disjoint from the other.
\end{quote}
\end{lem}
 
The next proposition establishes relations between the Borel rank
of the $G$-orbit of $x$ and the number $\gamma^G_{\star}(x)$.
The left inequality  is a direct consequence of Lemma \ref{char} 
and Proposition \ref{ef}, the right inequality follows from
Theorem \ref{orb} and Corollary \ref{complex}. 
 
\begin{prop} \label{ineq}
Let $x\in X$ and $\lambda$ be the multiplicative Borel rank of the orbit
$G x$ (i.e. $\lambda=min\{\zeta:\ Gx\in \Pi^0_{\zeta}\}$).
Then
$$ 
\gamma^G_{\star}(x)\le \lambda\le \gamma^G_{\star}(x)+3.
$$
In particular if $\lambda$ is a limit  ordinal, then
it is equal to $\gamma^G_{\star}(x)$. 
\end{prop}

The left inequality improves the analogous inequality obtained by Hjorth, 
who in fact proved that $\gamma^G_{\star}(x)\le \lambda +1$.
 
Proposition \ref{ineq} shows that the number $\gamma^G_{\star}(x)$ 
and the Borel rank of the orbit $G x$ can not differ very much.
Nevertheless we will show below that they can be different.
The corresponding example uses Lemma \ref{char}.
\bigskip
 
{\bf Example.}
Consider the conjugacy action of $S_{\infty}$ on itself.
It is shown in Theorem 1.8 of \cite{ivanov} that
any conjugacy class of $S_{\infty}$ belongs to $\Pi^0_3$
and there are conjugacy classes of Borel rank 3.
Let us prove that $\gamma_{\star}(f)=1$ for every
$f\in S_{\infty}$. 
Accordingly to Lemma \ref{char}, it suffices to show
that for any pair of conjugates $f$ and $g$ and any
finite set $c$ of natural numbers if
$V^G_c  f\cap V^G_c  g=\emptyset$
then we can separate $V^G_c  f$ from $V^G_c  g$
by an open or a closed subset of $S_{\infty}$
(in terms of that theorem $f=v^h$ and $g=v^{h_1}$ for
some $v,h,h_1 \in S_{\infty}$ with $h\in V^G_{\sigma}$
and $h_1 \in V^G_{\sigma_1}$).
We start with the following claim.
\parskip0pt 

{\em Claim.}
Let $f,g\in S_{\infty}$ be two conjugates and
$c\in [\omega]^{<\omega}$.
Then $V^G_c f$ and $V^G_c g$ are disjoint if and
only if there are $k\in c$ and $m\in \mathbb{Z}$ such that
$$
\Big(f^m(k)\in c \ \vee \ g^m(k)\in c\Big)\ \wedge \ f^m(k)\not=g^m(k).
$$

{\em Proof.}
It is well-known that $f$ and $g$ are conjugate if and only if
their cycle types are the same.
We have to consider only nonempty sets $c$.
Let $c=\{k_0,k_1,\ldots k_s\}$.
\parskip0pt

The proof of ($\Leftarrow$) is easy.
To prove the converse, assume that for all $k\in c$ and
$m\in {\mathbb Z}$, $\ f^m(k)\in c\ \vee\ g^m(k)\in c\ $ implies
$f^m(k)=g^m(k)$.
We are going to define some $h\in V^G_c$ so that $f^h=g$.
\parskip0pt

We proceed as follows.
For every $j\le s$ and every $m\in {\mathbb Z}$ we put
$$
h_0(f^m (k_j))=g^m (k_j) \ (\mbox{ in particular } h_0(k_j)=k_j).
$$
It follows from the assumptions that $h_0$ is a well-defined bijection
$$
h_0:\{f^m(k_j):\ j\le s,\ m\in {\mathbb Z}\}\to
\{g^m(k_j):\ j\le s,\ m\in{\mathbb Z}\}
$$
such that
$$
(f|_{\{f^m(k_j):\ j\le s,\ m\in {\mathbb Z}\}})^{h_0}=
g|_{\{g^m(k_j):\ j\le s,\ m\in {\mathbb Z}\}}.
$$
Now using the fact that $f$ and $g$ have the same cycle types,
we see that $h_0$ can be extended to a permutation
$h\in S_{\infty}$ so that $f^h=g$.
$\Box$
\bigskip

We can now finish the proof of the main statement.
By the claim we find $k,l\in c$ and $m\in \omega$
such that
$$
\Big(f^m(k)=l \ \vee \ g^m(k)=l\Big)\ \wedge \ f^m(k)\not=g^m(k).
$$
Without loss of generality we may assume that
$f^m(k)=l\ \wedge \ g^m(k)\not=l$.
Let $A_k^l$ be the set of all bijections
$$\sigma=\left(\begin{array}{l@{\ \ldots\ }l}
k\ \ a_1\ a_2\ & a_{m-2}\ a_{m-1}\\
a_1\   a_2\ a_3\ &  a_{m-1}\ l
\end{array}\right),$$
where $\{a_1,a_2,\ldots,a_{m-1}\}\subseteq \omega$.

Then the set $\hat{A}_k^l=\bigcup\{V^G_{\sigma}:\ \sigma\in A_k^l\}$ is
open, it contains $V^G_c f$ and is disjoint from $V^G_c g$.
$\Box$


\section{$\alpha$-Sets and admissible sets} 

Section 3.1 contains the main results of the paper. 
In Section 3.2 we give some straightforward construction of 
admissible sets which satisfy all our assumptions. 

\subsection{Main results} 

The main notions of this section 
(codability and constructibility in an admissible set) 
were defined in Section 1. 
The following theorem allows us to code $\alpha$-sets 
$B_{\alpha}(x,\sigma )$ in admissible sets.
Here we use terms which appear in Definitions \ref{D1} - \ref{F1} 
and Lemma \ref{constr}.  

\begin{thm}\label{con}
Let ${\mathbb{A}}$ be an admissible set.
Let $G<S_{\infty}$ be a closed subgroup and $X$ be
a Polish $G$-space with a basis $\{ A_i :i\in \omega\}$.
Suppose that $x$ is codable in ${\mathbb{A}}$ and the relation
$$ 
Imp(c,l,k)\quad \Leftrightarrow \big( c\in [\omega ]^{<\omega}\wedge
l,k\in \omega\wedge \ A_k\subseteq V^G_c  A_l\big)
$$ 
is $\Sigma$-definable in ${\mathbb{A}}$.
Then for every $\sigma\in S^G_{<\infty}$ and every countable 
$\alpha \in {\mathbb{A}}$ the set $B_{\alpha}(x,\sigma)$ 
is constructible in ${\mathbb{A}}$.
Moreover there is a $\Sigma$-definable in ${\mathbb{A}}$ binary 
function $u_x (\alpha ,\sigma )$ which finds 
a co-$\alpha$-multicode for $B_{\alpha }(x,\sigma )$.   
\end{thm}

{\em Proof.} 
We start with some preliminary remarks. 
First for every $c\in [\omega]^{<\omega}$
and $l\in \omega$ we define in ${\mathbb{A}}$ a function
$$
imp_{(c,l)}: \omega\to \{0,1\}\quad \mbox{ by  }\quad
imp_{(c,l)}(k)=1\ \mbox { iff }\ {\mathbb{A}}\models Imp(c,l,k). 
$$
Then  since $Imp$ is a $\Sigma$-relation in ${\mathbb{A}}$, we have
(by $\Sigma$-replacement): \parskip0pt

- $imp_{(c,l)}\in {\mathbb{A}}\ $,
${\mathbb{A}} \models B_{\Sigma}(1,imp_{(c,l)})\ $ and
$\ B_{imp_{(c,l)}}=V^G_c  A_l$; \parskip0pt

- $(0,imp_{(c,l)}))\in {\mathbb{A}}\ $,
${\mathbb{A}}\models B_{\Pi}(1,(0,imp_{(c,l)}))\ $ and
$\ B_{(0,imp_{(c,l)})}=X\setminus V^G_c A_l$. \\
As we noted in Section 1, the condition that $x$ is
$\Sigma$-codable in ${\mathbb{A}}$ implies that
$\langle [\omega]^{<\omega}, \subseteq\rangle$,
$\langle S^G_{<\infty}, \subseteq \rangle\ \in{\mathbb{A}}$.
We fix in ${\mathbb{A}}$ some bijective enumerations
$$
\rho :\omega\to [\omega ]^{<\omega}\quad
\quad \mbox{ and } \mu :\omega\to S^G_{<\infty}.
$$

Let $\sigma$ and $\alpha$ be as in the formulation
and let $d$, $c$ denote the domain and the range of
$\sigma$ respectively.
Let $F_1$ be a $\Sigma$-function which codes $x$ in
${\mathbb{A}}$. \parskip0pt

By the definition of $\alpha$-sets (see Definition \ref{ph}),
the set $B=B_{\alpha}(x,\sigma )$ can be naturally 
considered as an intersection of a pair of sets $C$ 
and $D$ (when $\alpha$ is limit we put $C$=$D$) of 
the form $C=\bigcap\limits_{i\in I}C_i$
and $D=\bigcap\limits_{i\in I}D_i$.
To define a multicode for $B$ we view $C$ and $D$ as
$$
C=X\setminus \bigcup\limits_{i\in I}(X \setminus C_i)
\quad\mbox{ and }
D=X\setminus \bigcup\limits_{i\in I}(X \setminus D_i) .
$$
Assuming that some co-multicodes for all $C_i$ and $D_i$,
$i\in\omega$, are already known, we will find in
${\mathbb{A}}$ appropriate multicodes $w$ and $v$ for
$\bigcup\limits_{i\in I}(X\setminus C_i)$ and
$\bigcup\limits_{i\in I}(X\setminus D_i)$ respectively.
It is worth noting that these multicodes will correspond
to ordinals appearing in Corollary \ref{complex} as levels
of Borel hierarchy.
Then the co-multicode $u=\bigwedge \big( (0,w),\ (0,v)\big)$
will correspond to $C\cap D$. \parskip0pt

Now we are ready to go into the details.
We define in ${\mathbb{A}}$ binary functions
$w_x (\alpha,\sigma)$ and $v_x (\alpha,\sigma)$ to the set
of multicodes and a function $u_x (\alpha,\sigma)$ 
\footnote{which in fact appears in the formulation of the theorem} 
to the set of co-multicodes by a formula depending 
on variables $\alpha ,\sigma ,u,v,w$ 
(realizing a simultaneous induction on ordinals $\alpha$) 
as follows:
$$
(\star)\
\big( w_x(\alpha,\sigma)=w\wedge v_x(\alpha,\sigma)=v \wedge
u_x(\alpha,\sigma)=u \big)
$$
$$
 \mbox{ iff }
$$
$$
\big( w, v \mbox{ are functions}\big) \
\wedge\  \big( u=\bigwedge ((0,w),\ (0,v))\big)\ \wedge
\bigg(\Theta_1\vee \Theta_2\vee \Theta_3 \bigg) ,
$$
where the formulas $\Theta_i$, for  $i=1,2,3$, are defined as follows.
\parskip0pt

$\Theta_1$ describes coding of $1$-sets:
$$
\Theta_1 =\ (\alpha =1)\wedge (dom[w]= dom[v]=\omega )
\wedge (\forall l\in \omega)\bigg( w(l) \mbox{ and } v(l)
\mbox{ are as follows:}
$$
$$
(l\in F_1(\sigma)\Rightarrow w(l)=(0,imp_{(c,l)}))
\wedge (l\not\in F_1(\sigma)\Rightarrow w(l)=mc_{\emptyset})\wedge 
$$
$$
(l\not \in F_1(\sigma)\Rightarrow v(l)=imp_{(c,l)})
\wedge (l\in F_1(\sigma)\Rightarrow v(l)=mc_{\emptyset})\bigg).
$$
Accordingly to this definition $w(l)$ is a co-$1$-multicode
of $X\setminus V^G_{c}\ A_l$ in the case
$V^G_{\sigma} x\cap A_l\not=\emptyset$.
Otherwise $w(l)$ is a $1$-multicode for $\emptyset$.
Eventually, $w$ is a $2$-multicode for the set
$\bigcup\{ X\setminus V^G_{c}\ A_l :V^G_{\sigma} x\cap A_l\not
=\emptyset\}$.
On the other hand $v$ is a $2$-multicode for the set
$\bigcup\{ V^G_{c}\ A_l :V^G_{\sigma}  x\cap A_l=\emptyset\}$.
\parskip0pt

The formula $\Theta_2$ tells us how to code
$\alpha$-sets at the successor step.
In this formula (see below) $w(l)$ is a co-multicode for the set
$X\setminus\bigcup\{ B_{\beta}(x,\sigma'):\sigma'\supseteq\sigma \wedge
dom[\sigma']=\rho (l)\}$
if $\ \rho (l)\supseteq dom[\sigma]\}$.
Otherwise $w(l)$ is a multicode for $\emptyset$.
Eventually, $w$ is a multicode for
$\bigcup\limits_{b\supseteq dom[\sigma]}
\bigg( X\setminus \bigcup\{B_{\beta}(x,\sigma'):\sigma'\supseteq\sigma\wedge
dom[\sigma']=b\}\bigg)$.\parskip0pt

On the other hand $v(l)$ is a co-multicode for
$X\setminus \bigcup\{ B_{\beta}(x,\sigma'):\sigma'\supseteq\sigma
\wedge rng[\sigma']=\rho (l)\}$
in the case $\rho (l)\supseteq rng[\sigma]$.
Otherwise $v(l)$ is a multicode for $\emptyset$.
Eventually, $v$ is a multicode for
$\bigcup\limits_{ a\supseteq rng[\sigma]}\bigg( X\setminus
\bigcup\{ B_{\beta}(x,\sigma'):\sigma'\supseteq\sigma\wedge rng[\sigma']
=a\}\bigg)$.

$$
\Theta_2 =\ (\exists \beta<\alpha)\bigg[ (\alpha=\beta +1)\wedge
(dom[w]=dom[v]=\omega ) \wedge
$$
$$
\wedge (\forall l\in \omega)\bigg( w(l)\mbox{ and }
v(l) \mbox{ are defined as follows: }
$$
$$
( \rho (l)\not\supseteq d\Rightarrow w(l)=mc_{\emptyset})\wedge
\big( \rho (l)\supseteq d\Rightarrow (\exists w'_l)( w(l)=(0,w'_l)\wedge
$$
$$
(w'_l\mbox{ is a function with }dom[w'_l ]=\omega )\wedge
(\forall j\in \omega )((\mu (j)\supseteq \sigma\wedge dom[\mu (j)]=\rho (l)
$$
$$
\Rightarrow w'_l (j)=u_x (\beta,\mu (j)) )\wedge (
(\mu (j)\not\supseteq \sigma\vee dom[\mu (j)]\not=\rho(l))\Rightarrow
w'_l(j)=mc_{\emptyset})))\big)\wedge
$$

$$
(\rho (l)\not\supseteq c\Rightarrow v(l)=mc_{\emptyset})\wedge
\big( \rho (l)\supseteq c\Rightarrow (\exists v'_l)( v(l)=(0,v'_l)\wedge
$$
$$
(v'_l\mbox{ is a function with } dom[v'_l ]=\omega )\wedge
(\forall j\in\omega )((\mu (j)\supseteq \sigma\wedge rng[\mu (j)]=\rho (l)
$$
$$
\Rightarrow v'_l (j)=u_x (\beta,\mu (j)) )\wedge (
(\mu (j)\not\supseteq \sigma \vee rng[\mu (j)]\not=\rho (l))\Rightarrow
v'_l (j)=mc_{\emptyset})))\big) \bigg)\ \bigg] .
$$

Finally, formula $\Theta_3$ settles the coding of $\alpha$-sets for
limit ordinals.

\begin{quote}
$\Theta_3=\ (\alpha$ is a limit ordinal)
$\wedge (w=v)\wedge (dom[w]=\alpha)\wedge (\forall\beta<\alpha)
\bigg( w(\beta)=(0,w'(\beta))$
where $w'(\beta)$ is a function defined on
$\omega$ such that $w'(\beta )(n)=u_x (\beta,\sigma)$
for every $n\in\omega\bigg)$.
\end{quote}
Thus $w(\beta )$, for every $\beta<\alpha$, is a co-multicode for
$X\setminus B_{\beta}(x,\sigma)$ and $w$ is a multicode for the union
$\bigcup\{X\setminus B_{\beta}(x,\sigma ):\beta<\alpha\}$.

Again we shall use the second recursion theorem
(Section 5.4 of \cite{barwise}) to see that $(\star)$ defines
a $\Sigma$-relation in ${\mathbb{A}}$.
Using induction and $\Sigma$-collection principle we conclude that
for every $\alpha \in {\mathbb{A}}$ and $\sigma\in S^G_{<\infty}$
the relation uniquely defines a co-multicode
$u_x (\alpha,\sigma)\in {{\mathbb{A}}}$ such that
$B_{\alpha}(x,\sigma)=B_{u_x (\alpha,\sigma)}$.
$\Box$
\bigskip

We have to prove the following technical lemma.
We apply the relation $\equiv$ defined in Section 1 (Definition \ref{eq}).

\begin{lem}\label{inj}
Let $G<S_{\infty}$ be a closed subgroup and $X$ be
a Polish $G$-space with a basis $\{ A_i :i\in \omega\}$.
Let ${\mathbb{A}}$ be an admissible set such that 
$Imp$ is $\Sigma$-definable in ${\mathbb{A}}$ and
$x, y$ are codable in ${\mathbb{A}}$.
Let $u_x$, $u_y$ denote the $\Sigma$-functions defined in the proof
of Theorem \ref{con}.
Then for every $\sigma, \delta\in S^G_{<\infty}$ with
$rng[\sigma]=rng[\delta]$
and every countable ordinal $\alpha\in {\mathbb A}$ we have
$$
B_{\alpha}(x,\sigma)=B_{\alpha}(y,\delta)
\quad \Rightarrow\quad
u_x(\alpha,\sigma)\equiv u_y(\alpha,\delta).
$$
\end{lem}

{\em Proof.}
We preserve the notation of the proof of Theorem \ref{con}.  
The proof is by induction on $\alpha$.
Assume that 
$$
(\star)\ rng[\sigma]=rng[\delta]\mbox{ and }
B_{\alpha}(x,\sigma)=B_{\alpha}(y,\delta).
$$
If $\alpha=1$ then  we apply the formula $\Theta_1$ from the proof 
of Theorem \ref{con}. 
Since $B_1(x,\sigma)=B_1(y,\delta)$, then for every $l\in\omega $,
$A_l$ intersects one of the sets $\ V_{\sigma} x, V_{\delta} y\ $
if and only if it also intersects the other.
Hence we have $w_x(1,\sigma)= w_y(1,\delta)$, $v_x(1,\sigma)= v_y(1,\delta)$.
Thus $u_x(1,\sigma)\equiv u_y(1,\delta)$.
\bigskip

For the sucessor step suppose that $\alpha=\beta +1$ 
and the implication 
$$
B_{\beta}(x,\sigma')=B_{\beta}(y,\delta')
\quad \Rightarrow\quad
u_x(\beta,\sigma')\equiv u_y(\beta,\delta').
$$
holds whenever $rng[\sigma']=rng[\delta']$.
We claim that 
$w_x(\alpha,\sigma)\equiv w_y(\alpha,\delta)$
and $v_x(\alpha,\sigma)\equiv v_y(\alpha,\delta)$. 
\parskip0pt

By Proposition \ref{=} the condition $(\star)$ implies 
that for every $a\supseteq c$ we have
$$
\quad \{B_{\beta}(x,\sigma'):\sigma'\supseteq\sigma, rng[\sigma']=a\}=
\{B_{\beta}(y,\delta'):\delta'\supseteq\delta, rng[\delta']=a\}. 
$$
Thus by the inductive assumption
for every $l\in \omega$ such that $\rho(l)\supseteq rng[\sigma]$ 
the sets 
$\{u_x(\beta,\sigma'):\sigma'\supseteq\sigma, rng[\sigma']=\rho(l)\}$ 
and $\{u_y(\beta,\delta'):\delta'\supseteq\delta, rng[\delta']=\rho(l)\}$ 
represent the same $\equiv$-classes. 
Hence $v_x(\alpha,\sigma)\ \equiv v_y(\alpha,\delta)$.
\parskip0pt 

On the other hand for every $b\supseteq dom[\sigma]$ there are 
$\sigma_1 \supset \sigma$ and $\delta_1 \supset\delta$
with $dom[\sigma_1 ]=b$ and $B_{\alpha}(x,\sigma_1 )=B_{\alpha}(y,\delta_1 )$.   
Then applying Lemma \ref{list}(c) (and the claim from its proof) 
we find $b_1\supseteq dom[\delta]$ (as $dom[\delta_1 ]$) such that 
$$
\quad \{B_{\beta}(x,\sigma'):\sigma'\supseteq\sigma, dom[\sigma']=b\}=
\{B_{\beta}(y,\delta'):\delta'\supseteq\delta, dom[\delta']=b_1\}. 
$$
By a similar argument we see that for every $b_1\supseteq dom[\delta]$ 
there is $b\supseteq dom[\sigma]$ such that
$$
\quad \{B_{\beta}(x,\sigma'):\sigma'\supseteq\sigma, dom[\sigma']=b\}=
\{B_{\beta}(y,\delta'):\delta'\supseteq\delta, dom[\delta']=b_1\}.
$$ 
Thus for every $l\in \omega$ such that $\rho(l)\supseteq dom[\sigma]$
we can find $l_1\in\omega$ (and vice-versa for every $l_1\in\omega$ 
with $\rho(l_1)\supseteq dom[\delta]$ there is $l\in\omega$ with
$\rho(l)\supseteq dom[\sigma]$)
such that the sets
$\{u_x (\beta,\sigma'):\sigma'\supseteq\sigma, dom[\sigma']=\rho(l)\}$ 
and $\{u_y(\beta,\delta'):\delta'\supseteq\delta, dom[\delta']=\rho(l_1)\}$
represent the same $\equiv$-classes. 
Therefore $w_x(\alpha,\sigma)\ \equiv w_y(\alpha,\delta)$.
By the formula $\Theta_2$ this finally yields
$u_x(\alpha,\sigma)\ \equiv u_y(\alpha,\delta)$.\bigskip

For the limit step suppose that $\alpha\in {\mathbb A}$ 
is a limit ordinal and  the implication
$$
B_{\beta}(x,\sigma)=B_{\beta}(y,\delta)
\quad \Rightarrow\quad
u_x(\beta,\sigma)\equiv u_y(\beta,\delta)
$$
holds for every $\beta<\alpha$.
Then $(\star)$ implies that for every $\beta<\alpha$
we have $B_{\beta}(x,\sigma)=B_{\beta}(y,\delta)$,
which by the inductive assumption gives
$u_x(\beta,\sigma)\equiv u_y(\beta,\delta)$.
Then we are done by the formula $\Theta_3$.
$\Box$
\bigskip

We shall now prove our main results
(which were formulated in Section 1).

\bigskip

{\em Proof of Theorem \ref{MaRe}}.
Let ${\mathbb{A}}$ be an admissible set such that
$\omega$ is realizable in it.
Let $G<S_{\infty}$ be a closed group, $X$ be a Polish $G$-space
with a basis $\{ A_i:i>0\}$ and $Imp$ be $\Sigma$-definable on ${\mathbb{A}}$.
We want to prove the following statements:
\begin{quote}
{\em (1) Let $x\in X$ be $\Sigma$-codable in ${\mathbb{A}}$.
Then for every $y\in X$, if $x,y$ are in the same
invariant Borel subsets of $X$ which are constructible
in ${\mathbb{A}}$ then for every $\alpha\le o({\mathbb{A}})$
they are in the same invariant $\Sigma_{\alpha}^0$-subsets of $X$.
\parskip0pt

(2) If $x,y$ are $\Sigma$-codable in ${\mathbb{A}}$ and they
belong to the same invariant Borel sets which are constructible
in ${\mathbb{A}}$ then they are in the same $G$-orbit.}
\end{quote}
Part (1) is a direct consequence of Theorem \ref{con}
and Lemma \ref{min}.\bigskip

(2) We shall use the back-and-forth arguments together
with $\Sigma$-reflection in ${\mathbb{A}}$.
We are going to construct a set $\Gamma$ of triples
$(n_i, \sigma_i, A_{i} )$ with the following properties
for every $i\in \omega$:

(a) $n_i\in \omega$, $\sigma_i\in S_{<\infty}^G$,
$rng[\sigma_{2i}]=n_{2i}$ and \parskip0pt

$dom[\sigma_{2i+1}]=n_{2i+1}$
(i.e. $rng[\sigma^{-1}_{2i+1}]=n_{2i+1}$);\parskip0pt

(b) $A_{2i}$ is a $V^G_{n_{2i}}$-invariant basic open
set containing $y$,\parskip0pt

$A_{2i+1}$ is a $V^G_{n_{2i+1}}$-invariant  basic open
set containing $x$; \parskip0pt

(c) $n_{i+1}>n_i$, $\sigma_{i+1}\supseteq \sigma_i$,
 $diam(A_{i+1})<2^{-(i+1)}$;
\parskip0pt

(d) $A_{2i+2}\subseteq A_{2i}$, and
$A_{2i+1}\subseteq A_{2i-1}$,\parskip0pt

$A_{2i+1}\subseteq V^G_{\sigma^{-1}_{2i}}  A_{2i}$ and
$A_{2(i+1)}\subseteq V^G_{\sigma_{2i+1}} A_{2i+1}$;\parskip0pt

(e)
$B_{\alpha}(x ,\sigma_{i})=B_{\alpha} (y, id_{rng[\sigma_{i}]})$

for every $\alpha < o({\mathbb{A}})$. \bigskip

We put $n_0 =0$, $\sigma_0 =\emptyset$ and let $A_0$ be
any $G$-invariant basic open set containing $y$. \parskip0pt

Suppose that we have already constructed all the
triples $(n_k, \sigma_k, A_{k})$, for $k\le 2i$.
In particular we have
$B_{\alpha}(x ,\sigma_{2i})=B_{\alpha} (y, id_{rng[\sigma_{2i}]})$
for every $\alpha\in o({\mathbb{A}} )$.
Then by Lemma \ref{list} (c) we also have
$B_{\alpha}(x ,id_{dom[\sigma_{2i}]})=B_{\alpha} (y, \sigma^{-1}_{2i})$,
for every $\alpha\in o({\mathbb{A}} )$.
\parskip0pt

Applying assumptions of the induction (in particular
$B_1 (x,\sigma_{2i})=B_1(y,id_{rng[\sigma_{2i}]})$)
we see that $V^G_{\sigma^{-1}_{2i}}  A_{2i}$ is
an $V^G_{dom[\sigma_{2i}]}$-invariant set containing $x$.
Let $A_{2i+1}\subseteq V^G_{\sigma^{-1}_{2i}}  A_{2i}$ be
any basic neighbourhood of $x$ such that $diam(A_{2i+1})<2^{-(2i+1)}$.
Then we define $n_{2i+1}$ to be any natural number
greater then $n_{2i}$ and covering $dom[\sigma_{2i}]$
such that $A_{2i+1}$ is $V^G_{n_{2i+1}}$-invariant.
\parskip0pt

We claim that there is some $\sigma_{2i+1}\in S^G_{<\infty}$
with $\sigma_{2i+1}\supseteq \sigma_{2i}$ and
$dom[\sigma_{2i+1}]=n_{2i+1}$ such that
$B_{\alpha} (x,id_{n_{2i+1}})=B_{\alpha} (y ,\sigma^{-1}_{2i+1})$
for every $\alpha < o({\mathbb{A}})$.
In other words we are looking for some
$\delta\in S^G_{<\infty}$ such that
$$(*)\
\delta\supseteq \sigma^{-1}_{2i},\
rng[\delta]=n_{2i+1}\mbox{ and }
B_{\alpha} (x ,id_{n_{2i+1}})=B_{\alpha} (y, \delta)\mbox{, for every }
\alpha\in o({\mathbb{A}} ).
$$
Suppose there is no $\delta$ satisfying $(*)$.
Then to every $\delta\in S^G_{<\infty}$,
satisfying $\delta\supseteq \sigma^{-1}_{2i}$ and
$rng[\delta]=n_{2i+1}$, we can assign some ordinal
$\beta_{\delta}\in {\mathbb{A}}$ so that
$B_{\beta_{\delta}} (x, id_{n_{2i+1}})\not= B_{\beta_{\delta}}(y, \delta)$.
By Lemma \ref{inj} the latter inequality is
equivalent to the relation
$u_x(\beta_{\delta},id_{n_{2i+1}})\not\equiv u_y(\beta_{\delta}, \delta)$.
By Definition \ref{eq} (and the discussion
after this definition) and Theorem \ref{con} this relation
can be expressed in ${\mathbb{A}}$ by a $\Sigma$-formula.
Since the set
$\{\delta\in S^G_{<\infty}:
\delta\supseteq \sigma^{-1}, rng[\delta]=n_{2i+1}\}$
is  an element of ${\mathbb{A}}$,
then by $\Sigma$-reflection in ${\mathbb{A}}$
(Section 1.4 of \cite{barwise}), there is an ordinal
$\beta\in {\mathbb{A}}$ such that
$u_x(\beta, id_{n_{2i+1}})\not\equiv u_y(\beta, \delta)$ for every
$\delta\supseteq \sigma^{-1}_{2i}$ with $rng[\delta]=n_{2i+1}$.
Therefore by the definition of the functions $u_x$ and $u_y$
(in the proof of Theorem \ref{con}) we have
$u_x(\beta+1, id_{dom[\sigma_{2i}]})\not\equiv u_y(\beta +1,\sigma^{-1}_{2i})$.
This by Lemma \ref{inj}  yields
$B_{\beta +1}(x, id_{dom[\sigma_{2i}]})\not=
B_{\beta +1}(y,\sigma^{-1}_{2i})$
which contradicts the assumptions. \parskip0pt

Then we take any $\delta$ satisfying $(\star)$ and
put $\sigma_{2i+1}=\delta^{-1}$.
At even steps we use the symmetric procedure.
\parskip0pt

Using the method just described we  define a
sequence $\{ \sigma_i :i<\omega\}$ of elements of $S^G_{<\infty}$.
Since $G$ is closed, by (b) and (c) there is $f\in G$ such that
$\bigcap\limits_i V^G_{\sigma_i}=\{ f\}$
and $fx\in \bigcap\limits_i V^G_{\sigma_{2i+1}}A_{2i+1}$.

Moreover by (a)-(d) we have
$V^G_{\sigma_{2i+1}} A_{2i+1}\subseteq
V^G_{\sigma_{2i}} A_{2i+1}\subseteq V^G_{n_2i}A_{2i}=A_{2i}$.
Therefore $\{ y\}=\bigcap\limits_i A_{2i}$
$=\bigcap \limits_i V^G_{\sigma_{2i+1}} A_{2i+1}=\{ fx\}$.
Thus $y=fx$.
$\Box$
\bigskip

This theorem has a corollary in the style of
Nadel's work \cite{nadel} (see also \cite{barwise}, Corollary 7.7.4).
It connects ranks considered in Section 2.2 with
$o({\mathbb{A}})$, the ordinal of ${\mathbb{A}}$.

\begin{prop}
Under the assumptions of this section (of Theorem \ref{con})
$\gamma^{G}_{\star}(x) \le  o({\mathbb A})$.
Moreover if $\gamma^{G}_{\star}(x) <  o({\mathbb{A}})$ then
$Gx$ is constructible in ${\mathbb A}$
(thus the Borel rank of $Gx$ is $<o({\mathbb{A}})$).
\end{prop}

{\em Proof.}
We apply the characterization of $\gamma^{G}_{\star}(x)$
from Lemma \ref{char}.
Take any $\sigma, \delta\in S^G_{<\infty}$ with
$rng[\sigma]=rng[\delta]=c$
such that $V^G_{\sigma}x\cap V^G_{\delta}x=\emptyset$.
It suffices to show that $V^G_{\sigma}x$ and $V^G_{\delta}x$ can
be separated by a Borel set constructible in ${\mathbb A}$
(its Borel rank is $< o({\mathbb{A}})$. \parskip0pt

First observe that if $x$ is codable in ${\mathbb A}$
with respect to $G$,
then it is also codable in ${\mathbb A}$ with respect
to any basic open subgroup $V_c^G$.
Indeed, let $F_1 :S_{<\infty}\to  {\mathbb A}$ be
a function coding $x$ in ${\mathbb A}$.
Since
$S^G_c=\{\sigma\in S^G_{<\infty}:(\forall n\in c\cap dom[\sigma ])(\sigma (n)=n)\}$
is an element of ${\mathbb A}$,
we can define in ${\mathbb{A}}$ a function coding $x$
with respect to $V^G_c$ by the following formula
$$
C_1(\sigma)=\left\{\begin{array}{l@{\quad\mbox{if}\quad}l}
\emptyset&\sigma\not\in S_c^G\\
F_1(\sigma)&\sigma\in S_c^G .
\end{array}
\right.
$$
It is clear that $C_1$ is an element of ${\mathbb A}$.\parskip0pt

In the following claim we consider elements of $G$ as
functions $\omega \rightarrow \omega$ with respect to
the chosen realization of $\omega$ in ${\mathbb{A}}$.
\parskip0pt

{\em Claim.} Let $x$ be codable in ${\mathbb A}$.
Then the set $\{g\in G:\ g\in {\mathbb A}\}$ is dense in $G$.
Moreover, if $g\in G$  is an element of ${\mathbb A}$,
then $gx$ is codable in ${\mathbb A}$.\parskip0pt

{\em Proof of Claim.} If $x$ is codable in ${\mathbb A}$ then
$S^G_{<\infty}$ is an element of ${\mathbb A}$.
We define on $S^G_{<\infty}$ a partial ordering $\le$ by the following
formula
$$
\sigma\le \delta\mbox{ iff } \sigma=\delta\vee
\Big( dom[\sigma]=dom[\delta]\ \wedge \
(\exists n\in dom[\sigma])
(\sigma(n)<\delta(n)\wedge 
$$
$$
(\forall k\in dom[\sigma ])(k<n \rightarrow\sigma(k)=\delta(k))\Big).
$$
Since this is a $\Delta_0$-formula, we see that
$\le$ is an element of ${\mathbb A}$.
Moreover for every $n$ the restriction of the ordering
$\ \le\ $ to the set
$\{\sigma\in S^G_{<\infty}:\ dom[\sigma]=n\}$
becomes a lexicographical well-order.
\parskip0pt

Now, take any $\sigma\in S^G_{<\infty}$.
We have to find $g\in G$ such that $\sigma\subseteq g$
and $g\in {\mathbb A}$.
We  define two increasing  sequences: $(k_n)$ of elements
of $\omega$ and $(\sigma_n)$ of elements of $S^G_{<\infty}$
by the following  scheme:
$$
\begin{array}l

\sigma_0=\sigma\\
k_0=min\{l:rng[\sigma_0]\subseteq l\}\\
\\

\sigma^{-1}_{2n+1}=
min_{\le}\{\delta\in S^G_{<\infty}:\ k_{2n}+2n+1=dom[\delta]
\ \wedge\ \sigma^{-1}_{2n}\subseteq \delta\}\\
k_{2n+1}=min\{l:dom[\sigma_{2n+1}]\subseteq l\}\\
\\
\sigma_{2n+2}=
min_{\le}\{\delta\in S^G_{<\infty}:\ k_{2n+1}+2n+2=dom[\delta]
\ \wedge\ \sigma_{2n+1}\subseteq \delta\}\\
k_{2n+2}=min\{l:rng[\sigma_{2n+2}]\subseteq l\}.
\end{array}
$$
We see that $g=\bigcup\limits_n \sigma_n$ is a permutation.
Since $G$ is closed, $g$ belongs to $G$.
On the other hand the definition of the functions $g$ and
$n\rightarrow\sigma_n$, $n\in \omega$, can be formalized by
a $\Sigma$-formula (by the second recursion theorem).
Since $\omega\in {\mathbb A}$, by $\Sigma$-replacement we have that
both the sequence $(\sigma_n)$ and $g$ are elements of ${\mathbb A}$.
\parskip0pt

Finally, if $g\in G$ and $g\in {\mathbb A}$, then
the operation $\ S^G_{<\infty}\ \to \  S^G_{<\infty}\ $:
$\ \sigma\to \sigma g\ $ defined by the $\Delta_0$-formula
$$
\sigma g\mbox{ is a finite partial function and }(\forall n\in dom[\sigma g])
(g(n)\in dom[\sigma]\ \wedge \ \sigma g(n)=\sigma((g(n)))
$$
is an element ${\mathbb A}$.\parskip0pt

Hence the function $F^g_1$ coding $gx$   can be defined
as an element of ${\mathbb A}$ by the following formula.
$$
F^g_1(\sigma)=\left\{\begin{array}{l@{\quad\mbox{if}\quad}l}
\emptyset&\sigma\not\in S^G_{<\infty}\\
F_1(\sigma g)&\sigma\in S^G_{<\infty}.
\end{array}
\right.
$$

Now we return to  the main statement.
By the claim there are $g\in V^G_{\sigma}$ and $f\in V^G_{\delta}$
such that $gx$ and $fx$  are codable in ${\mathbb A}$ with respect to $V_c$.
Then by Theorem \ref{MaRe}(2) (applied to $V^G_c$) we see that
$V^G_c gx=V^G_{\sigma}x$ and $V^G_c fx=V^G_{\delta}x$ can be
separated by Borel set constructible in ${\mathbb A}$.

To obtain the second part of the proposition note that
by Theorem \ref{orb} we see that $Gx$ is of the form
$B_{\alpha}(x,\emptyset )$ for some $\alpha \in{\mathbb{A}}$.
By Theorem \ref{con} this set is constructible in ${\mathbb{A}}$.
$\Box$
\bigskip

Theorem \ref{MaRe} suggests that in some situations we may
expect that $Gx$ is just the intersection of all
$G$-invariant Borel sets containing $x$ and constructible in ${\mathbb{A}}$.
We now show that under some additional assumption this
is really true.
This is the content of Theorem \ref{mor}:

\begin{thm}
Let $G$ be a closed subgroup of $S_{\infty}$,
$X$ be a Polish $G$-space, $t$ be a nice topology for $X$
and $\mathcal{B}$ be its nice basis.
Let $x\in X$ and let ${\mathbb{A}}$ be an admissible set such that
$x$ is codable in ${\mathbb{A}}$ with respect to $\mathcal{B}$.
Then the piece $C$ of the canonical partition
with respect to $\mathcal{B}$ with $x\in C$ coincides
with the orbit $G x$ if and only if $C$ is the intersection of
all invariant Borel sets containing $x$ and constructible in ${\mathbb{A}}$.
\end{thm}

{\em Proof.}
Since $C$ is invariant and Borel, the necessity is obvious.
To prove sufficiency we shall use  the notion of a {\em type}
introduced in \cite{basia}.

\begin{quote}
Let $H$ be an open subgroup of $G$ and $\hat{X}_0$
be an invariant $G_{\delta}$-subset of $X$ with respect to
the $t$-topology. \\
(1)A  family ${\cal F}\subseteq {\mathcal B}$ is called
an $H$-type in $\hat{X}_0$, if it is maximal with
respect to the following conditions:\parskip0pt

(a) $B$ is $H$-invariant, for any $B\in {\cal F}$;
\parskip0pt

(b) $\hat{X}_0\cap\bigcap {\cal F}\not=\emptyset$.\\
(2) An $H$-type ${\cal F}$ is called principal
if there is $B_{\cal F}\in {\cal F}$ such that
$B_{\cal F}\cap \hat{X}_0\subseteq B\cap \hat{X}_0$,
for every $B\in {\cal F}$.
We will say that $B_{\cal F}$ defines $H$.\\

\end{quote}

In paper \cite{basia} we prove the following characterization
of $G$-orbits in terms of types (Theorem 10).
\begin{quote}
Consider the canonical partition with respect to the
topology $t$.
A piece $\hat{X}_0$ of the canonical partition is
a $G$-orbit if and only if for any basic clopen subgroup $H<G$
any $H$-type of $\hat{X}_0$ is principal.
\end{quote}

Assuming (ii) we will show that every $H$-type of $C$
is principal for every clopen subgroup $H<G$.
Suppose the contrary.
Then there is some basic clopen $H<G$ and a non-principal $H$-type $\cal{F}$.
Then by Lemma 9 from \cite{basia}, the set
$D=\bigcap\limits_{g\in G}\bigg(g(\bigcup\limits_{B\in {\cal F}}(C\setminus B))\bigg)$
is  nonempty and invariant.
Since $C\cap \bigcap {\cal F}\not=\emptyset$, we have also $D\not=C$.
The Borel ranks of $D$ and $C\setminus D$
with respect to $\mathcal{B}$ are  $\le 4$
and one of these sets contains $x$.
This set by Lemma \ref{min} includes $B_4 (x,\emptyset)$.
Then we get a contradiction, since $B_4 (x,\emptyset)$
is constructible in ${\mathbb{A}}$ by Theorem \ref{con}.
$\Box$

\bigskip

In \cite{morozov} A.Morozov has proved the following theorem:
\begin{quote}
Let ${\mathbb{A}}$ be a {\em locally countable} admissible set
(i.e. $\omega <o({\mathbb{A}})$ and
${\mathbb{A}}\models (\forall s\not=\emptyset) (\exists f:\omega\rightarrow s)(f(\omega )=s)$).
Let $\phi$ be an $L_{\omega_1 \omega}$-sentence
for some language $L\in{\mathbb{A}}$.
Then $\phi$ is $\omega$-categorical if and only if
$\phi$ is complete with respect to all sentences
which belong to ${\mathbb{A}}$.
\end{quote}
Note that this theorem is quite similar to our Theorem \ref{mor}.
Indeed, consider the (logic) $S_{\infty}$-space $X_L$
of all $L$-structures.
Then identifying $L_{\omega_1 \omega}$-sentences $\psi$
with the corresponding $G_{\delta}$-set $\{ x: x\models\psi\}$
we see that sentences from ${\mathbb{A}}$ correspond
to $G_{\delta}$-sets constructible in ${\mathbb{A}}$.
Thus the condition that the set $C$ of models of $\phi$
cannot be divided by such $G_{\delta}$-sets means
that $C$ is an $S_{\infty}$-orbit. \parskip0pt

It is worth noting that our proof of Theorem \ref{mor}
is based on arguments which originally arose in model theory
(see \cite{barwise} and \cite{hodges}).

\subsection{Example of coding in admissible sets.}

Let $G$ be a closed subgoup of $S_{\infty}$, $(X,\tau )$
be a Polish $G$-space and ${\mathcal {A}}=\{A_l :l\in\omega\}$
be a  countable basis of $(X,\tau )$. \parskip0pt

To each $x\in X$ we assign an admissible set ${\mathbb{A}}_x$
such that $x$ is codable in ${\mathbb{A}}_x$. 
We start with the following two-sorted structure
$$
M_x=\langle \omega\cup S_{<\infty};
\ S^G_{<\infty}, Imp(\sigma ,k,l), Sat_x(\sigma, k)\rangle
$$
defined on the disjoint union of the set $\omega$ of natural numbers
and the set $S_{<\infty}$ of all bijections between
finite sets of natural numbers with:

\begin{enumerate}
\item The unary relation $S^G_{<\infty}$ for recognizing elements
of $S^G_{<\infty}$;
 
\item the ternary relation $Imp(\sigma ,k,l)$: 
$$ 
( \sigma\in S^G_{<\omega})\wedge
(l,k\in \omega )\wedge \ (A_k\subseteq V^G_{\sigma} A_l ); 
$$ 
\item the binary relation $Sat_x(\sigma,l)$ defined by
$Sat_x (\sigma,l)\Leftrightarrow  V^G_{\sigma} x\cap A_l\not=\emptyset$.
\end{enumerate}

\begin{prop}
Let ${\mathbb{A}}$ be an admissible set and $x\in X$.
The element $x$ is codable in  ${\mathbb{A}}$ so that 
$Imp (\sigma ,k,l)$ is definable for the corresponding 
realization of $(\omega ,<)$  
if and only if ${\mathbb{A}}$ is admissible above $M_x$ 
(i.e. $M_x \in {\mathbb{A}}$).
\end{prop}

{\em Proof.} 
$(\Rightarrow)$ According to Definition \ref{F1},
codability of  $x$  in ${\mathbb{A}}$
requires that ${\mathbb{A}}$ contains $\langle \omega, <\rangle$ or
its isomorphic copy.
Then as we have already noted in Section 1, 
${\mathbb{A}}$ also contains some copy of the structure 
$\langle S_{<\infty}, \subseteq \rangle$ 
. 
Let $F_1$ be the coding function for $x$. 
The predicate $S^G_{<\infty}$ is defined by the formula
$$
S^G_{<\infty}(\sigma)\ \mbox{ iff } \ F_1(\sigma)\not=\emptyset,
$$
so by $\Delta_0$-separation it also becomes an element of ${\mathbb{A}}$.
Finally, the relation
$Sat_x (\sigma, l)$ is also defined by the $\Delta_0$-fromula
$$
Sat_x (\sigma, l)\ \mbox{ iff } \  l\in F_1 (\sigma).
$$
\parskip0pt

The converse follows by similar arguments. \parskip0pt
$\Box$

\begin{cor}
Every $x\in X$ is codable in ${\mathbb{A}}_x ={\mathbf{Hyp}}(M_x)$.
\end{cor}

\bigskip

Institute of Mathematics, University of Wroc{\l}aw, \parskip0pt

pl.Grunwaldzki 2/4, 50-384 Wroc{\l}aw, Poland \parskip0pt

E-mail: biwanow@math.uni.wroc.pl

\end{document}